\documentclass[11pt]{article}
\usepackage{amssymb}
\usepackage{amscd}
\usepackage{times}
\usepackage[leqno]{amsmath}
\usepackage[all]{xy}

\usepackage{mathptmx}
\catcode`\@=11
\renewcommand\subsection{\@startsection{subsection}{2}{\z@}%
                                     {-3.25ex\@plus -1ex \@minus -.2ex}%
                                     {-0.01 mm}
                                     {\normalfont\large\bfseries}}

\catcode`\@=11
\renewcommand\subsubsection{\@startsection{subsubsection}{2}{\z@}%
                                     {-3.25ex\@plus -1ex \@minus -.2ex}%
                                     {-0.01 mm}
                                     {\normalfont\bfseries}}

\input epsf.tex

\newtheorem{Thm}{Theorem}[section]
\newtheorem{Lem}[Thm]{Lemma}
\newtheorem{Cor}[Thm]{Corollary}
\newtheorem{Prop}[Thm]{Proposition}
\newtheorem{Conj}[Thm]{Conjecture}
\newtheorem{example}[Thm]{Example}

\newtheorem{Def}[Thm]{Definition}

\newcommand{\E}{\mathcal{E}}

\newcommand{\Pro}{\mathbb{P}}
\newcommand{\N}{\mathbb{N}}
\newcommand{\C}{\mathbb{C}}
\newcommand{\R}{\mathcal{R}}

\newcommand{\md}{\operatorname{mod}}

\newcommand{\Hom}{\operatorname{Hom}} 
\newcommand{\Ext}{\operatorname{Ext}}

\newcommand{\bsm}{\begin{smallmatrix}}
\newcommand{\esm}{\end{smallmatrix}}

\newcommand{\orb}{\mathcal O}

\def\proof{\medskip\noindent {\it Proof --- \ }}

\def\resp{{\em resp. }}

\def\<{\langle\,}
\def\>{\,\rangle}

\def\eg{{\em e.g. }}
\def\SG{\mathfrak S}

\def\g{\mathfrak g}

\def\<{\langle}
\def\>{\rangle}
\def\n{{\mathfrak n}}
\def\le{\leqslant}

\def\ra{\rightarrow}

\def\1{\mathbf 1}
\def\Gr{{\rm Gr}}

\def\ii{{\mathbf i}}
\def\jj{{\mathbf j}}

\def\a{\alpha}

\def\L{\Lambda}
\def\l{\lambda}
\def\G{\Gamma}

\def\mod{{\rm mod}\,}

\def\Ext{{\rm Ext}}
\def\add{{\rm add}\,}

\def\Hom{{\rm Hom}}
\def\sHom{{\underline{\rm Hom}}}

\def\a{\alpha}
\def\b{\beta}

\def\Si{\Sigma}
\def\De{\Delta}

\def\x{{\mathbf x}}
\def\cqfd{\hfill $\Box$ \bigskip}

\def\vf{\varphi}
\def\sub{{\rm Sub\,}}
\def\fac{{\rm Fac\,}}

\def\pr{{\rm pr}}

\def\i{\iota}
\def\ra{\rightarrow}

\def\finex{\hfill$\diamond$}
\def\Om{\Omega}
\def\tx{\widetilde{x}}

\def\tf{\widetilde{f}}
\def\tg{\widetilde{g}}
\def\tB{\widetilde{B}}

\def\ie{{\em i.e.\ }}

\begin{document}

\title{\bf 
Partial flag varieties and preprojective algebras}
\author{C. Gei{\ss}, B. Leclerc and J. Schr\"oer}
\date{}
\maketitle

\begin{center} 
{\small \em Dedicated to Toshiaki Shoji on the occasion of his sixtieth birthday.}
\end{center}

\begin{abstract}
Let  $\L$ be a preprojective algebra of type $A, D, E$, and
let $G$ be the corresponding semisimple simply connected 
complex algebraic group.
We study rigid modules in subcategories $\sub Q$
for $Q$ an injective $\L$-module, and we introduce a mutation
operation between complete rigid modules
in $\sub Q$.
This yields cluster algebra structures on
the coordinate rings of the partial flag varieties attached to $G$.

\smallskip
\begin{center}
{\bf R\'esum\'e}
\end{center}

Soit  $\L$ une alg\`ebre pr\'eprojective de type $A, D, E$, et
soit $G$ le groupe alg\'ebrique complexe semi-simple et simplement
connexe correspondant.
Nous \'etudions les modules rigides des sous-cat\'egories
$\sub Q$ o\`u $Q$ d\'esigne un $\L$-module injectif, et nous introduisons
une op\'eration de mutation sur les modules rigides complets de $\sub Q$. 
Ceci conduit \`a des structures d'alg\`ebre amass\'ee sur les anneaux
de coordonn\'ees des vari\'et\'es de drapeaux partiels associ\'ees
\`a $G$.

\end{abstract}

\setcounter{tocdepth}{1}

{\small
\tableofcontents
}

\section{Introduction}


Let $G$ be a simple simply connected complex algebraic group
of simply-laced type $A$, $D$, $E$.
Let $P$ be a parabolic subgroup of $G$. 
The projective variety $P\backslash G$ is called a partial
flag variety. 
Let $\C[P\backslash G]$ be the multi-homogeneous coordinate
ring of $P\backslash G$ for its Pl\"ucker embedding into a product
of projective spaces. 

When $P=B$ is a Borel subgroup with unipotent radical $N$, 
$\C[B\backslash G]$ is equal to the coordinate ring of the 
base affine space $\C[N\backslash G]$. 
Let $w_0$ be  the longest element of the Weyl
group $W$ of $G$.
The double Bruhat cell $G^{e,w_0}$ \cite{FZ1} can be naturally 
identified to an open subset of $N\backslash G$.
Berenstein, Fomin and Zelevinsky \cite{BFZ} have shown
that $\C[G^{e,w_0}]$ has a nice cluster algebra structure, with 
initial seeds described explicitly in terms of the 
reduced expressions of $w_0$.

When $G=SL_n$ and $P$ is a maximal parabolic subgroup,
${\rm Gr}=P\backslash G$ is a Grassmann variety and Scott
\cite{S} has studied in detail a cluster
algebra structure on $\C[{\rm Gr}]$ 
(see also \cite{GSV}, where a cluster algebra on an open
Bruhat cell of ${\rm Gr}$ was previously introduced).

In \cite{GLS2}, we considered a cluster algebra structure
on $\C[N]$ coming from the one on $\C[G^{e,w_0}]$, and we
showed that it is deeply connected to the representation 
theory of the preprojective algebra $\L$ attached to the
Dynkin diagram of $G$.
We proved that every cluster of $\C[N]$ comes
in a natural way from a complete rigid $\L$-module
and that the cluster mutation lifts to 
a new operation on complete rigid modules 
that we also call mutation.
This allowed us to prove that every cluster monomial
is an element of the dual semicanonical basis, and
in particular that cluster monomials are linearly
independent. 

In this paper we consider an arbitrary
partial flag variety $P\backslash G$.
Let $N_P$ denote the unipotent radical of $P$.
We introduce a cluster algebra ${\cal A}\subseteq\C[N_P]$
whose initial seed is
described in terms of certain reduced expressions of $w_0$.
We then lift ${\cal A}$ to a cluster algebra
$\widetilde{\cal A}\subseteq\C[P\backslash G]$.
We conjecture that in fact 
${\cal A}=\C[N_P]$, and that
$\widetilde{\cal A}$ coincides with $\C[P\backslash G]$
up to localization with respect to certain generalized minors
(see \ref{subsecconj2} for a precise formulation).
We give a proof of the conjectures in type~$A$ and in type $D_4$.
We also give a complete classification of the algebras ${\cal A}$ 
(or equivalently $\widetilde{\cal A}$) which are of finite
type as cluster algebras. 
  
We arrived at the definition of ${\cal A}$ 
by studying certain subcategories of $\mod\L$.
Let $Q$ be an injective $\L$-module, and denote by 
$\sub Q$ the full subcategory of $\mod\L$ whose 
objects are the submodules of a finite direct sum of  
copies of~$Q$. 
This is a classical example of a homologically finite subcategory
(see \cite{AS,AR}).
It inherits from $\mod\L$ the structure of a Frobenius
category. 
We show that the relative syzygy functor of $\sub Q$
coincides with the inverse of the relative 
Auslander-Reiten translation.
This can be seen as a 2-Calabi-Yau property for $\sub Q$.
 
There is a natural correspondence between basic injective 
$\L$-modules $Q$ and conjugacy classes of parabolic subgroups $P$ of $G$. 
We show that the number of non-isomorphic indecomposable 
direct summands
of a rigid module in $\sub Q$ is at most 
the dimension of $P\backslash G$. 
A rigid module in $\sub Q$ with this maximal number of
summands is called complete. 
We prove that complete rigid modules exist in $\sub Q$ by
constructing explicit examples attached to certain
reduced words of $w_0$.
These modules give rise to the initial seeds from which
we define the cluster algebra ${\cal A}$.
If $X\oplus T\in\sub Q$ is a basic complete rigid module
with $X$ indecomposable and non-projective in $\sub Q$, 
we show that there exists a unique
indecomposable module $Y\in\sub Q$ non-isomorphic to $X$
such that $Y\oplus T$ is complete rigid. The module $Y\oplus T$ is
said to be obtained from $X\oplus T$ by a mutation.
We can attach to any complete rigid module an integer matrix  
encoding the mutations with respect to all its 
summands which are non-projective in $\sub Q$.
We show that these matrices follow the Fomin-Zelevinsky
matrix mutation rule.
This implies, as in \cite{GLS2}, that all cluster monomials
of ${\cal A}$ belong to the dual semicanonical basis.

To illustrate our results, we work out in detail in Section~\ref{sect12}
a simple but instructive example. 
We take for $G$ the group of type $D_n$ associated with 
a non-degenerate quadratic form on $\C^{2n}$.
The variety of isotropic lines is a smooth quadric
in $\Pro^{2n-1}(\C)$, which can be identified with $P\backslash G$
for a maximal parabolic subgroup $P$. 
We describe the corresponding subcategory $\sub Q$,
which is of finite type with $3n-4$ indecomposable modules,
and the cluster algebras ${\cal A} = \C[N_P]$ and
$\widetilde{\cal A} = \C[P\backslash G]$
which have finite cluster type $(A_1)^{n-2}$ in the Fomin-Zelevinsky
classification \cite{FZ2}.

In contrast, the cluster algebra structure on the coordinate 
ring of the Grassmann variety of isotropic subspaces of
dimension $n$ in $\C^{2n}$ is in general of infinite type. 
The only exceptions are $n=4$ and $n=5$, and they are described 
in detail in Section~\ref{sect13}.

We hope that, in general, our cluster algebra structure
$\widetilde{\cal A}$ on $\C[P\backslash G]$ will be helpful 
for studying total positivity and Poisson geometry on the
partial flag variety $P\backslash G$ (see \eg \cite{P,S,GSV}).

\bigskip
After this article was submitted for publication, the preprint
\cite{BIRS} of Buan, Iyama, Reiten and Scott appeared, which contains results similar to those
of Section~\ref{sect5} for more general subcategories. 
More precisely, it is proved in \cite[Chapter I, Th. 3.1]{BIRS}
that a theory of mutation for maximal rigid modules
can be developed for any extension closed functorially 
finite subcategory of $\md\,\L$.


\section{Coordinate algebras associated to parabolic subgroups}
\label{sectCoord}


In this section we fix our notation for algebraic groups
and flag varieties. 
We also recall some basic facts concerning coordinate
rings of partial flag varieties and their relation with
coordinate rings of unipotent radicals of parabolic subgroups.

\subsection{}
Let $\De$ be a Dynkin diagram of type $A,\,D,\,E$ 
with vertex set $I$.
We denote by $n$ the cardinality $|I|$ of $I$, and by $r$
the number of positive roots of $\Delta$.

Let $G$ be a simple simply connected complex algebraic group
with Dynkin diagram $\De$.
Let $H$ be a maximal  torus of $G$, and $B, B^-$ a pair of opposite 
Borel subgroups containing $H$ with unipotent radicals $N, N^-$. 

We denote by $x_i(t)\ (i\in I,\, t\in\C)$ the simple root subgroups
of $N$, and by $y_i(t)$ the corresponding simple root subgroups of
$N^-$.

The fundamental weights of $G$ are denoted by $\varpi_i\ (i\in I)$.
For a dominant integral weight $\l$ let $L(\l)$ be the
(finite-dimensional left) irreducible $G$-module with highest weight $\l$.
The $L(\varpi_i)$ are called fundamental representations.

Let $g \mapsto g^T$ be the involutive anti-automorphism of $G$ given by
\[
x_i(t)^T = y_i(t),\quad h^T=h,\qquad (i\in I, h\in H).
\] 
We denote by $L(\l)^*$ the right $G$-module obtained by twisting
the action of $G$ by this anti-auto\-mor\-phism.

The Weyl group of $G$ is denoted by $W$ and its longest element
by $w_0$. The Coxeter generators of $W$ are denoted by
$s_i \ (i\in I)$.
The length of $w\in W$ is denoted by $\ell(w)$.

The Chevalley generators of the Lie algebra $\g$ of $G$ are denoted
by $e_i, f_i, h_i\ (i\in I)$.
Here, the~$e_i$'s generate the Lie algebra $\n$ of $N$.

The coordinate ring $\C[N]$ is naturally endowed with a left action of $N$
\[
(x\cdot f)(n) = f(nx),\qquad (f\in\C[N],\ x,n\in N),
\]
and a right action of $N$
\[
(f\cdot x)(n) = f(xn),\qquad (f\in\C[N],\ x,n\in N).
\]
Differentiating these two actions we get left and right actions
of $\n$ on $\C[N]$. 
We prefer to write $e_i^\dag(f)$ instead of $f\cdot e_i$ for the
right action of the Chevalley generators.

\subsection{}
Throughout the paper we fix a non-empty subset $J$ of $I$ and
we denote its complement by $K=I\setminus J$.

Let $B_K$ be the standard parabolic 
subgroup\footnote{A more familiar notation for this parabolic subgroup would
be $P_K$, but we keep this symbol for denoting
certain projective modules over the 
preprojective algebra $\L$.} 
of $G$ generated by $B$
and the one-parameter subgroups 
\[
y_k(t),\qquad (k\in K,\,t\in\C).
\]
We denote by $N_K$ the unipotent radical of $B_K$.
In particular, we have $B_\emptyset=B$ and $N_\emptyset=N$.
On the other hand, when $K=I-\{j\}$ has $n-1$ elements,
$B_K$ is a maximal parabolic subgroup. 
It is known that every parabolic subgroup of $G$ is conjugate
to a unique standard parabolic subgroup~$B_K$ \cite[Cor. 11.17]{B}.

Let $B_K^-$ be the parabolic subgroup of $G$ generated by $B^-$
and the one-parameter subgroups 
\[
x_k(t),\qquad (k\in K,\,t\in\C).
\]
Then $N_K$ is an affine space of the same dimension as
the projective variety $B_K^-\backslash G$.
More precisely, the natural projection map
$$
\pi\colon G \ra  B_K^-\backslash G
$$ 
restricts to an embedding of $N_K$ into 
$B_K^-\backslash G$ as a dense open subset
(see \eg \cite[Prop. 14.21]{B}).

\begin{example}\label{exe1}
{\rm
Let $G=SL_5(\C)$, a group of type $A_4$. 
We choose for $B$ the subgroup of upper triangular 
matrices. 
Take $J=\{2\}$ and $K=\{1,3,4\}$.
Then $B_K^-$ and $N_K$ are the subgroups of $G$ 
with the following block form:  
\[
B_K^{-} =   
\begin{pmatrix} * & * & 0 & 0 & 0\cr
          * & * & 0 & 0 & 0 \cr
          * & * & * & * & * \cr
          * & * & * & * & * \cr
          * & * & * & * & * 
\end{pmatrix},
\qquad
N_K =
\begin{pmatrix} 1 & 0 & * & * & *\cr
          0 & 1 & * & * & * \cr
          0 & 0 & 1 & 0 & 0 \cr
          0 & 0 & 0 & 1 & 0 \cr
          0 & 0 & 0 & 0 & 1 
\end{pmatrix}. 
\] 
In this case $B_K^{-}$ is a maximal parabolic subgroup.
Let $(u_1,\ldots,u_5)$ be the canonical basis of $\C^5$.
We regard vectors of $\C^5$ as row vectors and
we let $G$ act on the {\em right}
on $\C^5$, so that the $k$th row of the matrix $g$
is $u_kg$.
Then $B_K^-$ is the stabilizer of the $2$-space
spanned by $u_1$ and $u_2$ for the induced 
transitive action of $G$ on the Grassmann
variety of $2$-planes of $\C^5$. Hence 
$B_K^-\backslash G$ is the Grassmannian $\Gr(2,5)$ 
of dimension~6. 
\finex
}
\end{example}

\subsection{}
The partial flag variety $B_K^-\backslash G$ can be 
naturally embedded as a closed subset in 
the product of projective spaces 
\[
\prod_{j\in J} \Pro(L(\varpi_j)^*)
\]
\cite[p.123]{LG}. 
This is called the Pl\"ucker embedding.
We denote by $\C[B_K^-\backslash G]$ the multi-homoge\-neous
coordinate ring of $B_K^-\backslash G$ coming from this 
embedding. 
Let $\Pi_J\cong \N^J$ denote the monoid of dominant integral 
weights of the form
\[
\l = \sum_{j\in J} a_j \varpi_j, \qquad (a_j\in \N).
\]
Then, $\C[B_K^-\backslash G]$ is a $\Pi_J$-graded ring with a natural 
$G$-module structure. 
The homogeneous component 
with multi-degree $\l\in\Pi_J$ is an irreducible $G$-module 
with highest weight $\l$.
In other words, we have
\[
\C[B_K^-\backslash G] = \bigoplus_{\l\in\Pi_J} L(\l).
\]
Moreover, $\C[B_K^-\backslash G]$ is generated by its subspace
$\bigoplus_{j\in J} L(\varpi_j)$.

In particular, $\C[B^-\backslash G]=\bigoplus_{\l\in\Pi} L(\l)$,
where the sum is over the monoid $\Pi$ of all dominant integral
weights of $G$. This is equal to the affine coordinate ring
$\C[N^-\backslash G]$ of the multi-cone $N^-\backslash G$  
over $B^-\backslash G$, that is, to the ring 
\[
\C[N^-\backslash G]=\{f\in\C[G]\mid f(n\,g)=f(g),\ n\in N^-,\ g\in G\}
\]
of polynomial functions on $G$ invariant under $N^-$.
The homogeneous component of degree $\l$ is the subspace
consisting of all functions $f$ such that 
$f(hg)=h^\l f(g)$ for $h\in H$ and $g\in G$.
We will identify $\C[B_K^-\backslash G]$ with the subalgebra
of $\C[N^-\backslash G]$ generated by the homogeneous elements
of degree $\varpi_j\ (j\in J)$.
\begin{example}\label{exe2}
{\rm We continue Example~\ref{exe1}. 
The Pl\"ucker embedding of the Grassmannian $\Gr(2,5)$
consists in mapping the $2$-plane $V$ of $\C^5$  
with basis $(v_1,v_2)$ to the line spanned by $v_1\wedge v_2$ 
in $\L^2\C^5$, which is isomorphic to the right  
$G$-module $L(\varpi_2)^*$ (remember that $G$ acts on the
right on $\C^5$).

This induces an embedding of $\Gr(2,5)$ into $\Pro(L(\varpi_2)^*)$.
The homogeneous coordinate ring for this embedding is 
isomorphic to the subring of $\C[G]$ generated by the
functions $g\mapsto \De_{ij}(g)$, where 
$\De_{ij}(g)$ denotes the
$2\times 2$ minor of $g$ taken on columns 
$i,j$ and on the first two rows.
The $\De_{ij}$ are called Pl\"ucker coordinates. 
As a $G$-module we have
\[
\C[\Gr(2,5)]=\bigoplus_{k\in \N} L(k\varpi_2),   
\]
where the degree $k$ homogeneous component $L(k\varpi_2)$
consists of the homogeneous polynomials of degree $k$
in the Pl\"ucker coordinates.
\finex}
\end{example} 

\subsection{}
Some distinguished elements of degree $\varpi_j$ in 
$\C[N^-\backslash G]$ are the generalized minors 
\[
\Delta_{\varpi_j, w(\varpi_j)},\qquad (w\in W),
\]
(see \cite[\S 1.4]{FZ1}).
The image of $N_K$ in $B_K^-\backslash G$ under the natural
projection is the open subset defined by the non-vanishing
of the minors $\De_{\varpi_j,\varpi_j}\ (j\in J)$.
Therefore the affine coordinate ring of $\C[N_K]$
can be identified with the subring of degree $0$ homogeneous 
elements in the localized ring 
$\C[B_K^-\backslash G][\De_{\varpi_j,\varpi_j}^{-1},\,j\in J]$.
Equivalently, $\C[N_K]$ can be identified with the quotient
of $\C[B_K^-\backslash G]$ by the ideal generated by the 
elements $\De_{\varpi_j,\varpi_j}-1\ (j\in J)$.

\begin{example}\label{exe3}
{\rm  We continue Example~\ref{exe1} and Example~\ref{exe2}. 
The coordinate ring of $\C[N_K]$ is isomorphic to the 
ring generated by the $\De_{ij}$ modulo the relation
$\De_{12}=1$. This description may seem  
unnecessarily complicated since 
$N_K$ is just an affine space of dimension $6$ and 
we choose a presentation with $9$ 
generators and the Pl\"ucker relations. 
But these generators are better adapted to the cluster
algebra structure that we shall introduce. 
\finex}
\end{example}

\subsection{}\label{ss93}
Let $\pr_J \colon \C[B_K^-\backslash G] \to \C[N_K]$ denote the 
projection obtained by modding out the ideal 
generated by the 
elements $\De_{\varpi_j,\varpi_j}-1\ (j\in J)$.
If $\C[B_K^-\backslash G]$ is identified as explained above
with a subalgebra of $\C[G]$, this map $\pr_J$ is nothing else
than restriction of functions from $G$ to $N_K$.
The restriction of $\pr_J$ to
each homogeneous piece $L(\l)\ (\l\in\Pi_J)$ of 
$\C[B_K^-\backslash G]$ is injective and gives an embedding of 
$L(\l)$ into $\C[N_K]$. 

We introduce a partial ordering on $\Pi_J$ by declaring that 
$\l\preccurlyeq \mu$ if and only if $\mu-\l$ is an $\N$-linear 
combination of weights
$\varpi_j\ (j\in J)$.
\begin{Lem}\label{Lem2.4}
Let $f\in\C[N_K]$.
There is a unique homogeneous element $\tf\in\C[B_K^-\backslash G]$
satisfying $\pr_J(\tf) = f$ and whose multi-degree is minimal
with respect to the above ordering $\preccurlyeq$. 
\end{Lem}
\proof
Let us first consider the algebra $\C[N]$.
Let $\l=\sum_i a_i\varpi_i$.
It is known that the subspace 
$\pr_I(L(\l))\subset \C[N]$ can be described as 
\[
\pr_I(L(\l))=\{f\in\C[N]\mid (e_i^\dag)^{a_i+1}f = 0,\ i\in I\}.
\]
In particular, $\C[N_K]$ can be identified with the subalgebra
\[
\{f\in\C[N]\mid  e_k^\dag f = 0,\ k\in K\}.
\]
For $\l=\sum_j a_j\varpi_j\in\Pi_J$, we then have
\[
\pr_J(L(\l)) = \{f\in\C[N_K]\mid  (e_j^\dag)^{a_j+1}f = 0,\ j\in J\}.
\]
Now given $f\in \C[N_K]$, we define 
$a_j(f) = \max\{s \mid  (e_j^\dag)^sf \not = 0\}$ and put
$\l(f)=\sum_{j\in J} a_j(f)\varpi_j$. 
Then $f\in\pr_J(L(\l(f)))$ and $\l(f)$ is minimal with this property.
Since the restriction of $\pr_J$ to $L(\l(f))$ is injective, we
see that there exists a unique $\tf$ as claimed :
this is the unique element of multi-degree
$\l(f)$ in $\C[B_K^-\backslash G]$ with $\pr_J(\tf)=f$.
\cqfd

For $f\in\C[N_K]$, let $a_j(f)$ and $\l(f)$ be as in the proof
of Lemma~\ref{Lem2.4}.
\begin{Lem}\label{lem25}
 Let $f,g\in\C[N_K]$. Then
$\widetilde{(f\cdot g)}=\tf\cdot\tg$.
If moreover, for every $j\in J$ we have 
$a_j(f+g)=\max\{a_j(f),\,a_j(g)\}$, 
then
\[
\widetilde{(f+g)} = \mu\tf + \nu\tg
\]
where $\mu$ and $\nu$ are monomials in the variables
$\De_{\varpi_j,\varpi_j}\ (j\in J)$ without common divisor.
\end{Lem}
\proof
The endomorphisms $e^\dag_j$ are derivations of $\C[N_K]$,
that is, we have
\[
e^\dag_j(fg)=e^\dag_j(f)g+fe^\dag_j(g),\qquad (f,g\in\C[N_K]).
\]
Therefore, by the Leibniz formula,
\[ 
(e^\dag_j)^{a_j(f)+a_j(g)}(fg)=(e^\dag_j)^{a_j(f)}(f)(e^\dag_j)^{a_j(g)}(g)\not = 0,
\]
and $(e^\dag_j)^{a_j(f)+a_j(g)+1}(fg)=0$. This proves the first
statement.
For the second statement, the additional assumption implies
that there exist unique monomials $\mu$ and $\nu$ in the 
$\De_{\varpi_j,\varpi_j}$ such that $\mu \tf$ and $\nu \tg$
have the same multi-degree as $\widetilde{(f+g)}$. 
Moreover $\mu$ and $\nu$ have no common factor.
The result follows.
\cqfd
 
\begin{example}\label{exe4}
{\rm We continue Examples~\ref{exe1}, \ref{exe2} and \ref{exe3}. 
For $1\le i<j\le 5$, let $D_{ij}=\pr_J(\De_{ij})$.
Thus, $D_{12}=1$, and for $(i,j)\not = (1,2)$, 
$\widetilde{D_{ij}}=\De_{ij}$.
It follows for example that, if $f=D_{13}D_{24}$
and $g=-D_{23}D_{14}$, then 
$\tf = \De_{13}\De_{24}$ and $\tg = -\De_{23}\De_{14}$.
Now we have the Pl\"ucker relation
\[
\De_{13}\De_{24}-\De_{23}\De_{14}=\De_{12}\De_{34}.
\]
This shows that $\widetilde{f+g}=\De_{34}$ is {\em not} of the
form $\mu\tf+\nu\tg$.
Here, we have $a_2(f)=a_2(g)=2$ but $a_2(f+g)=1$,
so that the assumption $a_2(f+g)=\max\{a_2(f),\,a_2(g)\}$
is not fulfilled.
\finex}
\end{example}


\section{The category $\sub Q_J$}


\subsection{}\label{subsec31}
Let $\L$ be the preprojective algebra of $\De$
(see \eg \cite{R,GLS0}).
This is a finite-dimensional selfinjective algebra over $\C$.
An important property of $\L$ is its $\Ext^1$-symmetry:
\[
\dim\Ext^1_\L(M,N)=\dim\Ext^1_\L(N,M),\qquad (M,N\in\mod\L).
\]
In particular, $\Ext^1_\L(M,N)=0$ if and only if 
$\Ext^1_\L(N,M)=0$. In this case, we say that the modules
$M$ and $N$ are {\em orthogonal}.
A $\L$-module $M$ is called {\em rigid} if it is orthogonal
to itself, that is, if $\Ext^1_\L(M,M)=0$.

The simple $\L$-module indexed by $i\in I$ is denoted by
$S_i$.
Let $P_i$ be the projective cover and $Q_i$ 
the injective
hull of $S_i$.
We denote by $\mu$ the Nakayama involution of $I$, defined
by $Q_i=P_{\mu(i)}$. 
(Equivalently, $\mu$ is characterized in terms of the
Weyl group $W$ by $s_{\mu(i)}=w_0s_iw_0$.) 

Let $\tau$ be the Auslander-Reiten translation of
$\mod \L$, and let $\Om$ be the syzygy functor. 
It is known that $\Om$ and $\tau^{-1}$ are isomorphic
as autoequivalences of the stable category
$\underline{\md}\,\L$ (see \eg \cite[\S 7.5]{GLS}). 

\subsection{}
Let 
$Q_J = \bigoplus_{j\in J} Q_j$ and
$P_J = \bigoplus_{j\in J} P_j$.
We denote by $\sub Q_J$ (\resp $\fac P_J$)
the full subcategory of $\mod \L$
whose objects are isomorphic to a submodule of a direct sum of
copies of $Q_J$ (\resp to a factor module of a direct sum of
copies of $P_J$). 
We are going to derive
cluster algebra structures on the coordinate algebras
$\C[N_K]$ and $\C[B_K^-\backslash G]$ 
by constructing a mutation operation on maximal
rigid modules in $\sub Q_J$.
We could alternatively use the subcategory $\fac P_J$,
but since this would lead to the same cluster structures
we shall only work with $\sub Q_J$.

For any unexplained terminology related to subcategories, we 
refer the reader to the introduction of \cite{AS}.
The next Proposition records some classical properties
of $\sub Q_J$ \cite[Prop. 6.1]{AS}. 
\begin{Prop}
The subcategory $\sub Q_J$ is closed under extensions, 
functorially finite, and has almost split sequences.
\cqfd
\end{Prop}

We shall denote by $\tau_J$ the {\em relative Auslander-Reiten translation}
of $\sub Q_J$.
Thus for $M$ an indecomposable object of $\sub Q_J$ which is 
not $\Ext$-projective, we have 
$N=\tau_J(M)$ if and only if there is an almost split sequence
in $\sub Q_J$ of the form
\[
0\to N\to E\to M\to 0.
\]
In this situation we also write $\tau_J^{-1}(N)=M$.

\subsection{}
For each $M$ in $\mod\L$, there is a unique submodule 
$\theta_J(M)$ minimal with respect to the property that
$M/\theta_J(M)$ is in $\sub Q_J$.
The natural projection $M\to M/\theta_J(M)$ is a minimal
left $\sub Q_J$-approximation of $M$ (see \cite[\S 4]{K}).
For $i\in I$, define 
$
L_i=P_i/\theta_J(P_i).
$

\begin{Prop}\label{prop2-1}
\begin{itemize}

\item[(i)]
The $L_i$ are the indecomposable $\Ext$-projective
and $\Ext$-injective modules in the subcategory $\sub Q_J$.

\item[(ii)]
The direct sum of the $L_i$ is a 
minimal finite cover and a minimal finite cocover of $\sub Q_J$.

\end{itemize}
\end{Prop}

\proof The $L_i$ are the indecomposable
$\Ext$-projectives by \cite[Prop. 6.3 (b)]{AS}. 
They are also $\Ext$-injective because of the $\Ext$-symmetry
of $\mod \L$. Since, by \cite[Prop. 6.3 (d)]{AS} the number of
indecomposable $\Ext$-injectives is the same as the number of
indecomposable $\Ext$-projectives, the $L_i$ are also all the 
indecomposable $\Ext$-injectives. This proves (i), and
(ii) follows from \cite[Prop. 3.1]{AS}.
\cqfd

\subsection{}
We introduce the {\em relative syzygy functor} $\Om_J$ of 
$\sub Q_J$. 
For $M\in\sub Q_J$, $\Om_J(M)$ is defined as the kernel 
of the projective cover of $M$ in $\sub Q_J$.

\begin{Lem}
For $M\in\sub Q_J$ we have 
\[
\Om_J(M)=\Om(M)/\theta_J(\Om(M)).
\]
\end{Lem}

\proof
Let $0 \to \Om(M) \stackrel{\iota}{\to} P \stackrel{\pi}{\to} M \to 0$ 
be the exact sequence defining $\Om(M)$. 
Since $M\in\sub Q_J$, we have 
$\theta_J(P)\subseteq \ker\pi=\iota(\Om(M)$,
hence $\theta_J(P)=\iota(N)$ for some $N\subseteq \Omega(M)$.
The induced sequence 
$0 \to \Om(M)/N \to P/\theta_J(P) \to M \to 0$
is exact, so $\Om(M)/N$ is isomorphic to a submodule of
$P/\theta_J(P)\in\sub Q_J$, therefore $\Om(M)/N\in\sub Q_J$
and $\theta_J(\Om(M))\subseteq N$.
Now, the induced sequence 
$0 \to \Om(M)/\theta_J(\Om(M)) \to P/\iota(\theta_J(\Om(M))) \to M \to 0$ 
is also exact, and since $\sub Q_J$ is closed under extensions,
this implies that $P/\iota(\theta_J(\Om(M)))\in \sub Q_J$,
hence $\theta_J(P)\subseteq\iota(\theta_J(\Om(M)))$,
that is, $N\subseteq\theta_J(\Om(M))$.
Therefore the sequence 
$0 \to \Om(M)/\theta_J(\Om(M)) \to P/\theta_J(P) \to M \to 0$
is exact in $\sub Q_J$, with middle term the projective cover of $M$ in
$\sub Q_J$, by Proposition~\ref{prop2-1}.
\cqfd

\begin{Prop}\label{prop2-2}
For any indecomposable non-projective module $M$ in $\sub Q_J$ 
we have
\[
\tau_J^{-1}(M)=\Om_J(M).
\]
\end{Prop}
\proof
Let $0\to M\to E \to \tau^{-1}(M) \to 0$ be an almost split sequence
in $\mod\L$.
Then by \cite[Corollary 3.5]{AS}, $\tau_J^{-1}(M)$ is a direct summand
of $\tau^{-1}(M)/\theta_J(\tau^{-1}(M))$. 
Now, it is known that $\tau^{-1}(M)/\theta_J(\tau^{-1}(M))$ is 
indecomposable \cite{BM,H}, see also \cite[Proposition 4.1]{K}.
Hence 
\[
\tau_J^{-1}(M)=\tau^{-1}(M)/\theta_J(\tau^{-1}(M))
=\Om(M)/\theta_J(\Om(M))=\Om_J(M).
\]
\cqfd

To summarize, the subcategory $\sub Q_J$ is a Frobenius category,
\ie it is an exact category with enough Ext-projectives, 
enough Ext-injectives, and moreover Ext-projectives and 
Ext-injectives coincide (see \cite{Ha}).
By Proposition~\ref{prop2-2}, the corresponding 
stable category $\underline{\sub} Q_J$ is a 2-Calabi-Yau
triangulated category (see \cite{Ke}).


\section{The category $(\sub Q_J)^\perp$}


A $\L$-module $M$ is said to be {\em orthogonal} to 
$\sub Q_J$ if $\Ext^1_\L(M,N)=0$ for every $N\in\sub Q_J$.
We now give a useful characterization of these modules.
 
\subsection{}
Let $\De_K$ be the subdiagram of $\De$ with
vertex set $K$.
This is in general a disconnected graph whose connected
components are smaller Dynkin diagrams.
Let $\L_K$ be the preprojective algebra of $\De_K$.
It is isomorphic to the direct product of the preprojective
algebras of the connected components of $\De_K$.
Although \cite{GLS2} only deals with preprojective
algebras associated with a single Dynkin diagram (the simple case),
it is easy to check that all the results extend to 
this more general ``semisimple case''.
In particular, maximal rigid $\L_K$-modules have $r_K$
non-isomorphic indecomposable direct summands, where $r_K$
denotes the number of positive roots of $\De_K$,
that is, the sum of the number of positive roots of each connected
component of $\De_K$.

In the sequel we shall always regard $\mod\L_K$ 
as the subcategory of $\mod\L$ whose objects
are the $\L$-modules supported on $K$.

\begin{Prop}\label{prop2-3}
As above let $K = I \setminus J$.
Then the following hold:
\begin{itemize}

\item[(i)]
A $\L$-module $M$ with no projective indecomposable direct summand satisfies 
$\Ext_\L^1(M,N) = 0$ for every $N\in\sub Q_J$ 
if and only if $M=\tau(U)$ with $U\in\mod\L_K$.
In other words, 
\[
(\sub Q_J)^\perp = \add(\tau(\mod\L_K)\cup\add\L).
\]

\item[(ii)]
A $\L$-module $N$ with no projective indecomposable direct summand satisfies 
$\Ext_\L^1(\tau(U),N) = 0$ for every $U\in\mod\L_K$ 
if and only if $N\in\sub Q_J$. 
In other words, 
\[
(\tau(\mod\L_K))^\perp = \add(\sub Q_J \cup \add\L).
\]
\end{itemize}
\end{Prop}

\proof 
(i)\quad By \cite[Prop 5.6]{AS}, $M$ satisfies the above property
if and only if 
\[
\Hom_\L(\tau^{-1}(M) , Q_J) = 0.
\] 
Now if $f\colon A \to Q_J$ is a nonzero homomorphism, then $f(A)$
is a submodule of $Q_J$ and its socle contains a module $S_j$
with $j\in J$. Therefore $A$ has a composition factor 
isomorphic to $S_j$ and $A\not \in \mod\L_K$.
Conversely, if $A$ has a composition factor
of the form $S_j$ with $j\in J$ we can find a submodule $B$
of $A$ such that the socle of $A/B$ contains only copies of 
$S_j$. Hence $A/B$ embeds in a sum of copies of $Q_j$, and
$\Hom_\L(A , Q_J) \not = 0$. 
So we have proved that $\Hom_\L(\tau^{-1}(M) , Q_J) = 0$ if and
only if all the composition factors of $U=\tau^{-1}(M)$
are of the form $S_k$ with $k\in K$, that is, $U\in\mod\L_K$.
This finishes the proof of $(i)$.

(ii)\quad
By part (i), if $N\in\sub Q_J$ then 
$\Ext_\L^1(\tau(U),N) = 0$ for every $U\in\mod\L_K$.
Conversely, if $N\in\mod\L$ is such that  
$\Ext_\L^1(\tau(U),N) = \Ext^1_\L(N,\tau(U)) = 0$ for every $U\in\mod\L_K$,
then by taking $U=S_k\ (k\in K)$ and using the classical
formula
\[
\Ext_\L^1(A,B)\simeq {\rm D}\sHom_\L(\tau^{-1}(B),A)
\]
we get that $\sHom_\L(S_k,N)=0$.
Therefore, every non-zero homomorphism from $S_k$ to $N$
factors through a projective-injective $\L$-module, which
can only be $Q_k$. This implies that $Q_k$ is a summand
of $N$, contrary to the hypothesis. Hence we have
that $\Hom_\L(S_k,N)=0$, and since this holds for every
$k\in K$, it follows that $N\in\sub Q_J$.
\cqfd

We can now give an alternative description of the modules $L_i$
of Proposition~\ref{prop2-1}. 
We denote by $q_k\ (k\in K)$ the indecomposable injective
$\L_K$-modules, that is, $q_k$ is the injective hull of
$S_k$ in $\mod\L_K$.

\begin{Prop} 
Let $i\in I$.
If $\mu(i)=j\in J$ then $L_i = Q_j=P_i$.
Otherwise, if $\mu(i)=k\in K$ then $L_i=\tau(q_k)$. 
\end{Prop}

\proof
Clearly, by definition of $L_i$, if $\mu(i)=j\in J$ then 
$L_i=Q_j=P_i$.
On the other hand, for every $U\in\mod\L_K$ we have 
\[
\Ext^1_\L(\tau(U),\tau(q_k)) \simeq 
\Ext^1_\L(U,q_k) =
\Ext^1_{\L_K}(U,q_k)=0, \qquad (k\in K).
\]
Hence, by Proposition~\ref{prop2-3} (ii), we have  
$\tau(q_k)\in\sub Q_J$.
Now, by Proposition~\ref{prop2-3} (i), $\tau(q_k)$ is $\Ext$-orthogonal
to every module $N\in\sub Q_J$, therefore $\tau(q_k)$, which is
indecomposable, is one of the Ext-projective modules $L_i$ of
$\sub Q_J$. Finally, since $\tau = \Omega^{-1}$, the head of
$\tau(q_k)$ is $S_i$ with $i=\mu(k)$, so $\tau(q_k)=L_i$.
\cqfd


\section{The functors $\E_i$ and $\E_i^\dag$}


In this section we introduce certain endo-functors $\E_i$ and
$\E_i^\dag$ of $\mod\L$
which we use to construct rigid $\L$-modules.
These functors should be seen as the lifts to $\mod\L$
of the maps $e_i^{\rm max}$ and $(e_i^\dag)^{\rm max}$
from $\C[N]$ to $\C[N]$ defined by
\[
e_i^{\rm max} f = (e_i^k/k!) f\quad \mbox{where} \quad k=\max\{j\mid e_i^j f\not =0\},
\]
and
\[
(e_i^\dag)^{\rm max} f = ({e_i^\dag}^l/l!) f\quad
\mbox{where}\quad l=\max\{j\mid {e_i^\dag}^j f\not =0\}.
\]
\subsection{}
For $M\in\mod\L$ and $i\in I$, let $m_i^\dag(M)$ denote the multiplicity 
of $S_i$ in the socle of $M$.
We define an endo-functor $\E^\dag_i$ of $\mod\L$ as follows.
For an object $M\in\mod\L$ we put 
\[
\E_i^\dag(M)=M/S_i^{\oplus m_i^\dag(M)}.
\]
So we get a short exact sequence
$$
0 \to S_i^{\oplus m_i^\dag(M)} \to M \to \E_i^\dag(M) \to 0.
$$
If $f : M \to N$ is a homomorphism, we can compose it
with the natural projection $N \to \E_i^\dag(N)$ to obtain 
$\widetilde{f} : M \to \E_i^\dag(N)$. The $S_i$-isotypic
component of the socle of $M$ is mapped to $0$ by $\widetilde{f}$, hence
$\widetilde{f}$ induces a homomorphism 
$\E_i^\dag(f) : \E_i^\dag(M) \to \E_i^\dag(N)$. Clearly, $\E_i^\dag$ is an additive functor.

\begin{Prop}\label{prop2-1n}
The functors $\E_i^\dag\ (i\in I)$ satisfy the following relations:
\begin{itemize}
\item[(i)]
$\E_i^\dag\E_i^\dag = \E_i^\dag$. 

\item[(ii)]  
$\E_i^\dag\E_j^\dag=\E_j^\dag\E_i^\dag$ if $i$ and $j$ are not connected
by an edge in $\De$.

\item[(iii)] 
$\E_i^\dag\E_j^\dag\E_i^\dag = \E_j^\dag\E_i^\dag\E_j^\dag$ if $i$ and $j$ are connected
by an edge in $\De$.
\end{itemize}
\end{Prop}

\proof Claim (i) is evident. 
Next consider the largest submodule $M_{i,j}$ of $M$ whose 
composition factors are all isomorphic to $S_i$ or $S_j$.
If $i$ and $j$ are not connected by an edge in $\De$ then 
$\Ext^1_\L(S_i,S_j)=\Ext^1_\L(S_j,S_i)=0$,
hence $M_{i,j}$ is semisimple, and one clearly has
\[
\E_i^\dag\E_j^\dag(M)=\E_j^\dag\E_i^\dag(M)=M/M_{i,j}.
\]
If $i$ and $j$ are connected by an edge in $\Delta$,
an elementary calculation in type $A_2$ shows that
\[
\E_i^\dag\E_j^\dag\E_i^\dag(M_{i,j}) =
\E_j^\dag\E_i^\dag\E_j^\dag(M_{i,j}) = 0.
\]
It follows that
\[ 
\E_i^\dag\E_j^\dag\E_i^\dag(M) = \E_j^\dag\E_i^\dag\E_j^\dag(M) = M/M_{i,j}.
\]
\cqfd

Let $(i_1,\ldots,i_k)$ be a reduced word for $w\in W$.
Proposition~\ref{prop2-1n} implies that 
the functor $\E_{i_1}^\dag\cdots \E_{i_k}^\dag$
does not depend on the choice of the reduced word.
We shall denote it by $\E_w^\dag$.

\subsection{}
For $M\in\mod\L$ let $m_i(M)$ be the multiplicity
of $S_i$ in the top of $M$.
We define an endo-functor $\E_i$ of $\mod\L$ as follows.
For an object $M\in\mod\L$ we define $\E_i(M)$
as the kernel of the surjection 
$$
M \to S_i^{\oplus m_i(M)}.
$$
If $f\colon M \to N$ is a homomorphism,  
$f(\E_i(M))$ is contained in $\E_i(N)$, and
we define 
$$
\E_i(f)\colon \E_i(M) \to \E_i(N)
$$ 
as the restriction of $f$ to
$\E_i(M)$.
Clearly, $\E_i$ is a functor.
Alternatively, we could define $\E_i$ as the composition
of functors 
\begin{equation}\label{Sdual}
\E_i = S\E_i^\dag S,
\end{equation}
where $S$ is the self-duality of $\mod\L$ introduced in \cite[\S 1.7]{GLS}.
This shows immediately that the functors $\E_i$ satisfy 
analogues of Proposition~\ref{prop2-1n}. 
In particular, we have for every $w\in W$ a well-defined functor
$\E_w$.

\begin{Prop}\label{prop2-2n}
If $M$ is rigid then $\E_i(M)$ and $\E_i^\dag(M)$ are rigid.
\end{Prop}
\proof 
Let $M$ be a rigid module with dimension vector $\b$.
We use the geometric characterization of rigid 
$\L$-modules given in \cite[Cor. 3.15]{GLS2}, namely,
$M$ is rigid if and only if its orbit $\orb_M$ in the
variety $\L_\b$ of $\L$-modules of dimension vector $\b$
is open.

Let $p=m_i(M)$ and $\b'=\b-p\a_i$. 
In  \cite[\S 12.2]{L91}, Lusztig considers, for $q\in\N$,
the subvariety $\L_{\b,i,q}$ of $\L_\b$ consisting of all $\L$-modules
$N$ with $m_i(N)=q$. He introduces the variety $Y$ 
of triples $(t,s,r)$ where $t\in \L_{\b',i,0}$,
$s\in\L_{\b,i,p}$ and $r: t \to s$ is an injective homomorphism
of $\L$-modules. Let $p_1 : Y\to \L_{\b',i,0}$ and 
$p_2: Y\to \L_{\b,i,p}$ denote the natural projections. 
It is easy to check that 
$p_1p_2^{-1}(\orb_N)=\orb_{\E_i(N)}$, 
the orbit of $\E_i(N)$ in $\L_{\b'}$.
By \cite[Lemma 12.5]{L91}, $p_2$ is a principal $G_{\b'}$-bundle,
and the map $(t,s,r) \mapsto (t,r)$ is a locally trivial 
fibration $p'_1 : Y \to \L_{\b',i,0} \times J_0$
with fibre isomorphic to $\C^m$ for a certain $m\in\N$.
Here, $J_0$ stands for the set of all injective graded
linear maps from $V_{\b'}$ to $V_\b$.
It follows that if $\orb_M$ is open in $\L_\b$ then
$p'_1p_2^{-1}(\orb_M)$ is open in $\L_{\b',i,0} \times J_0$,
and $\orb_{\E_i(M)}$ is open in $\L_{\b',i,0}$. 
Since $\L_{\b',i,0}$ has the same pure dimension
as $\L_{\b'}$ \cite[Th. 12.3]{L91}, 
we conclude that $\orb_{\E_i(M)}$ is open in $\L_{\b'}$.
This proves the claim for $\E_i(M)$. 
The result for  $\E_i^\dag(M)$ then follows from (\ref{Sdual}).
\cqfd

If $M$ is a rigid $\L$-module with dimension vector $\b$
then $\overline{\orb_M}$ is an irreducible component
of $\L_\b$.
So the functors $\E_i$ yield maps which associate to
an irreducible component containing a rigid module, other 
irreducible components, also containing a rigid module.
By \cite{KS}, these maps correspond to certain arrows
of Lusztig's coloured graph for the enveloping algebra $U(\n)$
of Dynkin type $\De$ (see \cite[\S 14.4.7, \S 14.4.15]{L93})
or, equivalently, to certain paths in Kashiwara's crystal
graph of $U(\n)$.

\subsection{}
The functors $\E_i$ and $\E_i^\dag$ can be used to construct rigid 
$\L$-modules. 
Start with $Q=\oplus_{i\in I}Q_i$.
This is a projective-injective module, hence a rigid module. 
Let $\ii=(i_1,\ldots,i_r)$ be a reduced word for $w_0$. 
We define a rigid module $T_\ii^\dag$ by induction
as follows. 
Put ${}^\dag T^{(r+1)}=Q$, and 
\[
{}^\dag T^{(k)} = \E_{i_k}^\dag\left({}^\dag T^{(k+1)}\right) \oplus  Q_{i_k},\qquad (1\le k\le r).
\]
Then define $T_\ii^\dag := {}^\dag T^{(1)}$.
\begin{example}\label{ex2-3}
{\rm
In type $A_3$, let $\ii=(2,1,3,2,1,3)$. We have
\[
Q_1=\def\objectstyle{\scriptstyle} \xymatrix@-1.2pc 
{&&3\ar[ld]\\&2\ar[ld]&\\1&&},\quad
Q_2=\def\objectstyle{\scriptstyle} \xymatrix@-1.2pc 
{&\ar[ld]2\ar[rd]&\\1\ar[rd]&&\ar[ld]3\\&2&},\quad
Q_3=\def\objectstyle{\scriptstyle} \xymatrix@-1.2pc 
{1\ar[rd]&&\\&2\ar[rd]&\\&&3}.
\]
Here the graphs display the structure of the injective $\L$-modules.
Thus, $Q_1$ is a uniserial module with socle $S_1$, top $S_3$
and middle layer $S_2$.
(In the sequel we shall frequently use graphs of this type to
represent examples of $\L$-modules.)
Now,
\[
{}^\dag T^{(6)}=\E_3^\dag(Q) \oplus Q_3=\E_3^\dag(Q_3)\oplus Q =
\def\objectstyle{\scriptstyle} \xymatrix@-1.2pc {1\ar[rd]&\\&2}\ \oplus Q.
\]
Next,
\[
{}^\dag T^{(5)}=\E_1^\dag({}^\dag T^{(6)}) \oplus Q_1 
= 
\def\objectstyle{\scriptstyle} \xymatrix@-1.2pc {1\ar[rd]&\\&2}\ 
\oplus\ 
\def\objectstyle{\scriptstyle} \xymatrix@-1.2pc {&\ar[ld]3\\2&}  
\oplus Q.
\]
Similarly we get
\[
{}^\dag T^{(4)}= 
\def\objectstyle{\scriptstyle}\xymatrix@-1.2pc{1}\ 
\oplus\  
\def\objectstyle{\scriptstyle}\xymatrix@-1.2pc{3}\ 
\oplus\ 
\def\objectstyle{\scriptstyle} \xymatrix@-1.2pc {&\ar[ld]2\ar[rd]&\\1&&3} 
\oplus\ Q,
\qquad
{}^\dag T^{(3)}= 
\def\objectstyle{\scriptstyle}\xymatrix@-1.2pc{1}\ 
\oplus\ 
\def\objectstyle{\scriptstyle} \xymatrix@-1.2pc {&\ar[ld]2\\1&} 
\oplus\ 
\def\objectstyle{\scriptstyle} \xymatrix@-1.2pc {1\ar[rd]&\\&2} 
\oplus\ Q,
\]
\[
{}^\dag T^{(2)}=
\def\objectstyle{\scriptstyle} \xymatrix@-1.2pc {1\ar[rd]&\\&2} 
\oplus\ 
\def\objectstyle{\scriptstyle}\xymatrix@-1.2pc{2}\ 
\oplus\ 
\def\objectstyle{\scriptstyle} \xymatrix@-1.2pc {&\ar[ld]3\\2&} 
\oplus\ Q,
\qquad
{}^\dag T^{(1)}=
\def\objectstyle{\scriptstyle}\xymatrix@-1.2pc{1}\ 
\oplus\  
\def\objectstyle{\scriptstyle}\xymatrix@-1.2pc{3}\ 
\oplus\ 
\def\objectstyle{\scriptstyle} \xymatrix@-1.2pc {&\ar[ld]2\ar[rd]&\\1&&3} 
\oplus\ Q = T_\ii^\dag.
\]
\finex}
\end{example}

\subsection{}\label{ss24}
Similarly we can define a rigid module $T_\ii$ inductively
as follows. 
Put $T^{(r+1)}=Q$, and 
\[
T^{(k)} = \E_{i_k}\left( T^{(k+1)}\right) \oplus  P_{i_k},
\qquad (1\le k\le r).
\]
Then define 
$$
T_\ii :=  T^{(1)}.
$$
In fact, we have $T_\ii=ST_\ii^\dag$.
In particular, $T_\ii$ and $T_\ii^\dag$ have the same number
of indecomposable direct summands, and their endomorphism algebras
have opposite Gabriel quivers.

\subsection{}
Let $W_K$ be the subgroup of $W$ generated by the $s_k\ (k\in K)$.
Let $w_0^K$ denote the longest element of $W_K$. 
We have $\ell(w_0^K)=r_K$.

\begin{Prop}\label{lem4-4}
For any $\L$-module $M$, we have
\[
\E_{w_0^K}^\dag(M)=M/\theta_J(M).
\]
\end{Prop}

\proof
It is easily checked that for every $k\in K$,
$\E_{w_0^K}^\dag(q_k) = 0$.
This implies that for every $N\in\mod\L_K$ we have
$\E_{w_0^K}^\dag(N) = 0$.
Therefore, if $N$ is the largest submodule of $M$
contained in $\mod\L_K$, we have 
$\E_{w_0^K}^\dag(M) = M/N$,
and the result follows.
\cqfd

Proposition~\ref{lem4-4} shows that 
the functor $\E_{w_0^K}^\dag :\mod\L \to \sub Q_J$  
is left adjoint to the inclusion functor
$\iota_J : \sub Q_J \to \mod\L$ (see \eg \cite[\S2, p.17-18]{AR}).
In other words, there are isomorphisms
\[
\Hom_\L(M,N) \stackrel{\sim}{\to} \Hom_{\sub Q_J}(\E_{w_0^K}^\dag(M),N)
\]
functorial in $M\in\mod\L$ and $N\in\sub Q_J$.

In particular, taking $K=\{k\}$, we see that
$\E_k^\dag$ is left adjoint to the inclusion functor
$\iota_{I-\{k\}}$.
Similarly, $\E_k$ is right adjoint to
the inclusion functor of the subcategory 
$\fac P_{I-\{k\}}$.


\section{Generalized minors and maximal rigid modules}


\subsection{}
For $i\in I$ and $u,v\in W$, let $\De_{u(\varpi_i),v(\varpi_i)}$
denote the generalized minor introduced by Fomin and 
Zelevinsky \cite[\S 1.4]{FZ1}.
This is a regular function on $G$. 
We shall mainly work with the restriction of this 
function to $N$, that we shall denote by $D_{u(\varpi_i),v(\varpi_i)}$.
It is easy to see that $D_{u(\varpi_i),v(\varpi_i)}=0$
unless $u(\varpi_i)$ is less or equal to $v(\varpi_i)$ in the Bruhat 
order, and that $D_{u(\varpi_i),u(\varpi_i)}=1$ for
every $i\in I$ and $u\in W$ (see \cite[\S 2.3]{BFZ}).

\subsection{}
In \cite[\S 9]{GLS2} we have associated to every 
$M\in\mod\L$ a regular function $\varphi_M\in\C[N]$,
encoding the Euler characteristics of the varieties
of composition series of $M$.
In particular, one has 
\begin{equation}\label{equa1}
\varphi_{Q_i}=D_{\varpi_i,w_0(\varpi_i)},\qquad (i\in I).
\end{equation}
More generally, it follows from \cite[Lemma 5.4]{GLS}
that for $u, v\in W$ we have
\begin{equation}\label{equa2}
\varphi_{\E^\dag_u\E_v Q_i} = D_{u(\varpi_i),vw_0(\varpi_i)}.
\end{equation}

\subsection{}\label{ss25}
To a reduced word $\ii = (i_1,\ldots,i_r)$ of $w_0$,
Fomin and Zelevinsky have attached a sequence of minors
$\De(k,\ii)$ defined as follows (see \cite[\S 5.3]{GLS}). 
Here $k$ varies over $[-n,-1]\cup [1,r]$, 
\[
\De(-j,\ii) = \De_{\varpi_j,w_0(\varpi_j)},\qquad
(j\in[1,n]),
\]
and 
\[
\De(k,\ii) = \De_{\varpi_{i_k},v_{>k}(\varpi_{i_k})},\qquad
(k\in[1,r]),
\]
where $v_{>k}=s_{i_r}s_{i_{r-1}}\cdots s_{i_{k+1}}$.
We shall denote by $D(k,\ii)$ the restriction of
$\De(k,\ii)$ to $N$.
There are $n$ indices $l_j\ (1\le j\le n)$ in $[1,r]$ such that 
$D(l_j,\ii)=D_{\varpi_j,\varpi_j}=1$.
The remaining $r-n$ indices in $[1,r]$ are called $\ii$-exchangeable.
The subset of $[1,r]$ consisting of these $\ii$-exchangeable indices
is denoted by $e(\ii)$. 
The functions $D(k,\ii)\ (k\in[-n,-1]\cup e(\ii))$ form one of the
initial clusters for the cluster algebra $\C[N]$ \cite{BFZ,GLS}.

\begin{Prop}\label{prop2-4n}
Let $\ii = (i_1,\ldots,i_r)$ be a reduced word for $w_0$.
Then the following hold:
\begin{itemize}

\item[(i)]
$T_\ii$ is a basic complete rigid $\L$-module.

\item[(ii)] 
We can denote the $r$ indecomposable direct summands of 
$T_\ii$ by 
$T_k$ in such a way that 
\[
\varphi_{T_k}=D(k,\mu(\ii)),\qquad (k\in[-n,-1]\cup e(\mu(\ii))),
\]
where $\mu(\ii) = (\mu(i_1),\ldots,\mu(i_r))$ is the image of $\ii$
under the Nakayama involution (see \S\ref{subsec31}).

\end{itemize}
\end{Prop}

\proof
One has $\E_i(P_j)=P_j$ if $j\not = i$.
Moreover, since $P_i$ has a simple socle, all its
submodules are indecomposable.
It then follows from the inductive definition of $T_\ii$ 
that the indecomposable direct summands of $T_\ii$ are the
projective modules $P_i$ and all the nonzero modules
of the form
\[
M_k=\E_{i_1}\E_{i_2}\cdots\E_{i_k}(P_{i_k}),\qquad (1\le k\le r).
\]
Using Equation~(\ref{equa2}) and the
fact that $P_i=Q_{\mu(i)}$, we thus get immediately 
\[
\varphi_{M_k} = 
D_{\varpi_{\mu(i_k)},u_{\le k}w_0(\varpi_{\mu(i_k)})},
\] 
where $u_{\le k} = s_{i_1}\cdots s_{i_k}$.
Now,
\[
u_{\le k}w_0=w_0 s_{i_r}\cdots s_{i_{k+1}} w_0
=s_{\mu(i_r)}\cdots s_{\mu(i_{k+1})},
\]
and it follows that $\varphi_{M_k} = D(k,\mu(\ii))$.
Note that the $n$ indices $k\in[1,r]$ which are not 
$\mu(\ii)$-exchan\-gea\-ble 
give $\varphi_{M_k} = 1$,
hence the corresponding modules $M_k$ are zero.
Therefore, taking into account the $n$ modules $Q_i$, we
see that $T_\ii$ has exactly $r$ pairwise non-isomorphic 
indecomposable direct summands
and projects under the map $M\mapsto \varphi_M$ to the
product of cluster variables 
\[
D(k,\mu(\ii)),\qquad (k\in[-n,-1]\cup e(\mu(\ii))).
\] 
\cqfd

Proposition~\ref{prop2-4n} gives a representation-theoretic construction
of all the initial clusters of $\C[N]$.
Note that this description is different from that of \cite{GLS}
which was only valid for the reduced words $\ii$
adapted to an orientation of the Dynkin diagram.
It is known (see \cite[Remark 2.14]{BFZ})
that all these clusters are related to each other
by mutation.
Therefore, using \cite{GLS} and 
\cite{GLS2} we know that the exchange matrix of
the cluster 
$$
\{D(k,\mu(\ii)) \mid k\in [-n,-1]\cup e(\mu(\ii))\}
$$ 
coincides with the exchange matrix of the rigid module $T_\ii$.
Hence, using \cite{BFZ}, we have an easy combinatorial
rule to determine the exchange matrix of $T_\ii$ for
every $\ii$, or in other words the Gabriel quiver of its
endomorphism algebra.

\subsection{}\label{ss26}
For $u, v\in W$ and $i\in I$ we have
\[
S(\E^\dag_u\E_vQ_i)= \E_u\E^\dag_vP_i=\E^\dag_v\E_uQ_{\mu(i)}.
\]
It follows that 
\[
\varphi_{S(\E^\dag_u\E_vQ_i)}=D_{v(\varpi_{\mu(i)}),uw_0(\varpi_{\mu(i)})}.
\]
This, together with Proposition~\ref{prop2-4n}, implies that
\begin{Cor}\label{Cor62}
The rigid module $T^\dag_\ii$ is basic and complete in $\mod \L$.
Its $r-n$ non-projective indecomposable direct summands $T^\dag_k=ST_k$
satisfy
\[
\varphi_{T^\dag_k}=D_{u_{\le k}(\varpi_{i_k}), w_0(\varpi_{i_k})},
\qquad (k\in e(\ii)),
\]
where $u_{\le k} = s_{i_1}\cdots s_{i_k}$.
\cqfd
\end{Cor}

\begin{example}{\rm
We continue Example~\ref{ex2-3}.
We have 
\[
D_{u_{\le 6}(\varpi_3),w_0(\varpi_3)}=1,\qquad 
D_{u_{\le 5}(\varpi_1),w_0(\varpi_1)}=1,\qquad
D_{u_{\le 4}(\varpi_2),w_0(\varpi_2)}=1,
\]
\[ 
D_{u_{\le 3}(\varpi_3),w_0(\varpi_3)}=D_{s_2s_1s_3(\varpi_3),w_0(\varpi_3)}
=\varphi_{S_1},
\]
\[
D_{u_{\le 2}(\varpi_1),w_0(\varpi_1)}=D_{s_2s_1(\varpi_1),w_0(\varpi_1)}
=\varphi_{S_3},
\]
\[
D_{u_{\le 1}(\varpi_2),w_0(\varpi_2)}=D_{s_2(\varpi_2),w_0(\varpi_2)}
=\varphi_{\def\objectstyle{\scriptscriptstyle} 
\xymatrix@-1.2pc {&\ar[ld]2\ar[rd]&\\1&&3}}.
\]
}
\finex
\end{example}


\section{Complete rigid modules in $\sub Q_J$}


\subsection{}
For a $\L$-module $T$, we denote by $\Si(T)$ the number of 
non-isomorphic indecomposable direct summands of $T$.

\begin{Prop}
Let $T$ be a rigid module in $\sub Q_J$. 
We have
$\Si(T)\le r-r_K$.
\end{Prop}

\proof
We may assume without loss of generality that $T$ is basic.
Let $T'$ be a basic rigid $\L_K$-module without projective summands.
Then by \cite{GS} we know that $\Si(T')\le r_K-|K|$, and by \cite{GLS}
we can assume that this upper bound is achieved, that is, 
we can assume that $T'$ has $r_K-|K|$ indecomposable summands.
It follows that $\tau(T')$ is a basic rigid $\L$-module with
$r_K-|K|$ indecomposable summands.
By Proposition~\ref{prop2-3}, no summand of $\tau(T')$ 
belongs to $\sub Q_J$, since such a summand would have to be
the $\tau$-translate of a projective $\L_K$-module. 
Therefore the $\L$-module
$T''=T\oplus\tau(T')\oplus Q_K$
is basic and rigid.
Here we use that for $k\in K$, the injective module $Q_k$ is
not in $\sub Q_J$.
By \cite{GS} we have $\Si(T'')\le r$, hence 
\[
\Si(T) + (r_K - |K|) + |K| = \Si(T) + r_K \le r.
\]
\cqfd 

Note that $r-r_K= \dim N_K = \dim B_K^-\backslash G$.

\begin{Def}{\rm
A rigid module $T$ in $\sub Q_J$ is called {\em complete}
if $\Si(T)=r-r_K$.
}
\end{Def} 

\subsection{}\label{ss52}
We can write $w_0=w_0^Kv_K$ with 
$\ell(w_0^K)+\ell(v_K)=\ell(w_0)$.
Therefore there exist reduced words $\ii$ for $w_0$ 
starting with a factor $(i_1,\ldots,i_{r_K})$ 
which is a reduced word for $w_0^K$.

\begin{Prop}\label{prop5-4}
Let $\ii = (i_1,\ldots,i_r)$ be a reduced word for $w_0$ such that 
$(i_1,\ldots,i_{r_K})$ is a reduced word for $w_0^K$.
Then $T^\dag_\ii$ has $r-r_K$ non-isomorphic indecomposable direct summands in 
$\sub Q_J$.
\end{Prop}

\proof
The indecomposable direct summands of $T^\dag_\ii$ are the 
$Q_i\ (i\in I)$ and the modules
\[
M_j^\dag=\E^\dag_{i_1}\E^\dag_{i_2}\cdots\E^\dag_{i_j}(Q_{i_j}),\qquad (j\in e(\ii)).
\]
By Corollary~\ref{Cor62}, these $r$ summands are pairwise non-isomorphic.
For $j>r_K$, by Proposition~\ref{lem4-4} we have $M^\dag_j \in \sub Q_J$.
Taking into account the $n$ non-exchangeable elements, we have
$r-n-r_K$ such modules.
Moreover, again by Proposition~\ref{lem4-4}, we have
\[
N^\dag_k:=\E^\dag_{i_1}\E^\dag_{i_2}\cdots\E^\dag_{i_{r_K}}(Q_k)
=\E^\dag_{w_0^K}Q_k\in \sub Q_J,
\qquad (k\in K).
\]
Let $s_k=\max\{s\le r_K \mid i_s = k\}$. 
Then $M^\dag_{s_k}=N^\dag_k$ belongs to $\sub Q_J$.
Finally, since for $j\in J$ we have obviously $Q_j\in\sub Q_J$,
we have found $r-n-r_K+|K|+|J|=r-r_K$ summands of $T^\dag_\ii$ in
$\sub Q_J$.
\cqfd

Proposition~\ref{prop5-4} shows that there exist complete
rigid modules in $\sub Q_J$, and gives a recipe to construct
some of them.
Note that the modules $L_i\ (i \in I)$ 
occur as direct summands 
of every complete rigid module $T$ in $\sub Q_J$.
This follows immediately from the fact that the $L_i$ are 
Ext-projective and Ext-injective in $\sub Q_J$.

\subsection{}
Proposition~\ref{prop5-4} is based on properties of 
the functor $\E^\dag_{w_0^K}$. 
A related result is given by the next proposition.
Following Iyama, we call a $\L$-module $T$ 
{\em maximal 1-orthogonal in} $\sub Q_J$ if
$$
\add T = \{ M \in \sub Q_J \mid \Ext^1_\L(T,M)=0 \}.
$$
Clearly a maximal 1-orthogonal module $T$ is {\em maximal rigid in} $\sub Q_J$,
that is, it is rigid, and if $X \in \sub Q_J$ is such that
$T \oplus X$ is rigid then $X \in \add T$.

\begin{Prop}\label{1orth}
Let $T$ be a maximal rigid module in $\mod \L$.
Then $\E^\dag_{w_0^K}(T)$ is maximal 1-ortho\-go\-nal in $\sub Q_J$.
\end{Prop}

\proof
By Proposition~\ref{prop2-2n}, $\E^\dag_{w_0^K}(T)$ is rigid, so if
$M \in \add \E^\dag_{w_0^K}(T)$ then 
$$
\Ext^1_\L(\E^\dag_{w_0^K}(T),M) = 0.
$$
Conversely, let us assume that $M \in \sub Q_J$ is such that
$\Ext^1_\L(\E^\dag_{w_0^K}(T),M) = 0$.

Let $f\colon T' \to M$ be a minimal right $\add T$-approximation.
Since $\L \in \add T$, $f$ is surjective and putting $T'' = \ker f$
we get an exact sequence 
$$
0\to T'' \to T' \to M \to 0.
$$
By Wakamatsu's lemma, it
follows that $\Ext^1_\L(T,T'')=0$ (see \cite[Lemma 3.1]{AR}).
By \cite[Th. 6.4]{GLS2}, we know that $T$ is maximal 1-orthogonal
in $\mod\L$, hence $T'' \in \add T$.

Now, since $M \in \sub Q_J$ we have $\theta_J(M)=0$
and $\theta_J(T'') \simeq \theta_J(T')$.
Therefore we get an exact sequence
\[
0\to\E^\dag_{w_0^K}(T'') \to \E^\dag_{w_0^K}(T') \to M \to 0.
\]
But $\E^\dag_{w_0^K}(T'') \in \add \E^\dag_{w_0^K}(T)$, so by assumption
$$
\Ext^1_\L(\E^\dag_{w_0^K}(T''),M)=\Ext^1_\L(M,\E^\dag_{w_0^K}(T''))=0
$$ 
and the sequence splits.
Hence $M$ is a direct summand of $\E^\dag_{w_0^K}(T')$ and
$M \in\add \E^\dag_{w_0^K}(T)$, as required.
\cqfd

\subsection{}
It follows from \cite[Cor. 8.7]{GLS4} that every 
maximal rigid module in $\sub Q_J$ is complete.
Proposition~\ref{1orth} then implies that for every maximal rigid
module $T\in\mod\L$, the module $\E^\dag_{w_0^K}(T)$ is 
complete rigid in $\sub Q_J$.
We will not make use of this result in the sequel.


\section{Mutations in $\sub Q_J$}\label{sect5}


\subsection{}
Let $X\oplus T$ be a basic complete rigid module in $\sub Q_J$
with $X$ indecomposable and not Ext-projective.
Let $V$ be a basic rigid module in $\mod\L_K$ without projective
summand, and let us assume that $\Si(V)=r_K-|K|$.
Then, by Proposition~\ref{prop2-3}, 
$X\oplus T\oplus \tau(V)\oplus Q_K$
is a basic maximal rigid module in $\mod\L$.

Let $f : T' {\to} X$ be a minimal right $\add(T)$-approximation.
Since $T$ contains as a summand
the module $L=\oplus_j L_j$, which is a cover of $\sub Q_J$,
the map $f$ is surjective.
Let $Y$ be the kernel of $f$, so that we have a short exact
sequence
\begin{equation}\label{ses1}
0 \to Y \stackrel{g}{\to} T' \stackrel{f}{\to} X \to 0
\end{equation}
where $g$ denotes inclusion.
Since $T'$ belongs to $\sub Q_J$ which is closed under taking
submodules, $Y$ also belongs to $\sub Q_J$.
Using the dual of \cite[Proposition 5.6]{GLS2}, 
we get that $g$ is a minimal left $\add(T)$-approximation,
$T \oplus Y$ is basic rigid, $Y$ is indecomposable and
not isomorphic to $X$.
Again by Proposition~\ref{prop2-3}, we have that 
$Y\oplus T\oplus \tau(V)\oplus Q_K$
is a basic maximal rigid module in $\mod\L$.
By \cite{GLS2} this is the mutation of 
$X\oplus T\oplus \tau(V) \oplus Q_K$
in the direction of $X$.
This shows in particular that $f$ and $g$ are in fact
$\add(T\oplus \tau(V) \oplus Q_K)$-approximations,
and that $Y$ is the unique indecomposable module in $\sub Q_J$
non-isomorphic to $X$ and such that $T\oplus Y$ is complete rigid.
By \cite{GLS2} we have another short exact sequence
\begin{equation}\label{ses2}
0 \to X \stackrel{i}{\to} T'' \stackrel{h}{\to} Y \to 0
\end{equation}
where $h$ and $i$ are 
$\add(T\oplus \tau(V) \oplus Q_K)$-approximations.
Since $\sub Q_J$ is closed under extensions, $T''$
is in $\sub Q_J$.
Hence  $h$ and $i$ are also $\add(T)$-approximations,
and in particular $T''$ does not depend on the choice of $V$.
It also follows from \cite[Corollary 6.5]{GLS2} that 
$\dim \Ext^1_\L(X,Y) = 1$, hence the short exact sequences
(\ref{ses1}), (\ref{ses2}) are the unique non-split short
exact sequences between $X$ and $Y$.
To summarize, we have obtained the following
\begin{Prop}
Let $X\oplus T$ be a basic complete rigid module in $\sub Q_J$
with $X$ indecomposable and not $\Ext$-projective.
There exists a unique indecomposable module $Y\not\simeq X$ in 
$\sub Q_J$ such that $Y\oplus T$ is a basic complete rigid 
module in $\sub Q_J$.
Moreover, $\dim \Ext^1_\L(X,Y) = 1$ and if (\ref{ses1}), (\ref{ses2}) 
are the unique non-split short exact sequences between $X$ and $Y$,
then $f, g, h, i$ are minimal $\add(T)$-approximations.
\cqfd
\end{Prop}

In the situation of the above Proposition, we say that 
$Y\oplus T$ is the {\em mutation of $X\oplus T$ in the
direction of $X$}, and we write
$\mu_X(X\oplus T)=Y\oplus T$.
Let $U=X\oplus T$.
Let $B=B(U\oplus \tau(V)\oplus Q_K)^\circ$ be the exchange
matrix of the maximal rigid module 
$U\oplus \tau(V)\oplus Q_K$. 
The $(r-r_K)\times (r-r_K-n)$ submatrix of 
$B$ whose rows are labelled by the indecomposable summands 
of $U$ and columns by the indecomposable 
summands of $U \setminus L$ is called 
{\em the exchange matrix of $U$}, and denoted
$B(U)^\circ$.
By the discussion above, it contains all the information 
to calculate the mutations of $U$ in all the $r-r_K-n$
directions.
Moreover, the submatrix of $B$ whose rows are labelled by 
the summands of $\tau(V)\oplus Q_K$ and columns
by the indecomposable 
summands of $U \setminus L$ has all its entries
equal to $0$.
Taking into account the results of \cite{GLS2},
this implies the

\begin{Prop}
Let $U$ be a basic complete rigid module in $\sub Q_J$.
Write $U=U_1\oplus\cdots\oplus U_{r-r_K}$, where the last $n$ summands
are $\Ext$-projective. 
Let $k\le r-r_K-n$. 
We have 
\[
B(\mu_{U_k}(U))^\circ = \mu_k(B(U)^\circ),
\]
where in the right-hand side $\mu_k$ denotes the 
Fomin-Zelevinsky matrix mutation.
\cqfd
\end{Prop}


\section{Cluster algebra structure on $\C[N_K]$}\label{sect10}


\subsection{}
The next proposition relates the coordinate ring $\C[N_K]$
to the subcategory $\sub Q_J$.

\begin{Prop}\label{prop7-2}
The algebra $\C[N_K]$ is isomorphic to the subspace of $\C[N]$ spanned by 
the set
\[
\{ \vf_M \mid M\in\sub Q_J\}.
\]
\end{Prop}
\proof
By \ref{ss93}, 
$\C[N_K]$ is isomorphic to the (non direct) sum 
of the subspaces
\[
\pr_I(L(\l)),\qquad (\l\in\Pi_J),
\]
of $\C[N]$.
It was shown in \cite[Th. 3]{GLS1} that, 
if $\l=\sum_{j\in J}a_j\varpi_j$
and $Q(\l)=\oplus_{j\in J}Q_j^{\oplus a_j}$, then 
$\pr_I(L(\l))$ is spanned by the $\varphi_M$ where $M$
runs over all submodules of $Q(\l)$.
The result follows. 
\cqfd

\subsection{}\label{sect9.2}
Following \ref{ss52}, we introduce
the set $R(w_0,K)$ of reduced words $\ii=(\ii',\ii'')$
for $w_0$ starting with a reduced word $\ii'$ of $w_0^K$.
By Proposition~\ref{prop5-4}, if $\ii\in R(w_0,K)$,
the basic complete rigid module $T_\ii$
has a unique direct summand which is a 
complete rigid module in $\sub Q_J$.
We shall denote this summand by $U_\ii$.
If $\jj=(\jj',\jj'')$ is another word in $R(w_0,K)$, 
we can pass from $\ii'$ to $\jj'$ by a sequence of $2$-moves
and $3$-moves, and similarly from $\ii''$ to $\jj''$
by a sequence of $2$-moves and $3$-moves.
Using \cite[Remark 2.14]{BFZ}, this implies that 
$T_\ii$ and $T_{(\jj',\ii'')}$ are connected by a sequence
of mutations. Moreover, the definition of the modules
$T_\ii$ shows that all these mutations will leave unchanged
the direct summand $U_\ii\in\sub Q_J$.
Similarly, $T_{(\jj',\ii'')}$ and $T_\jj$ are connected
by a sequence of mutations and it is plain that they all take 
place in $\sub Q_J$. Hence the modules $U_\ii$ and $U_\jj$
are connected by a sequence of mutations in $\sub Q_J$.
To summarize, we have proved that
\begin{Lem}
The set $\R_J$ of basic complete rigid modules in $\sub Q_J$ which can
be reached from $U_\ii$ by a finite sequence of mutations does not
depend on the choice of\, $\ii\in R(w_0,K)$.\cqfd
\end{Lem}
Now, exactly as in \cite{GLS2}, we can project $\R_J$
on $\C[N_K]$ using the map $M\mapsto \varphi_M$.
More precisely, put $d_K=r-r_K=\dim N_K$, and for
$T=T_1\oplus \cdots \oplus T_{d_K}\in\R_J$, 
let $\x(T)=(\varphi_{T_1},\ldots,\varphi_{T_{d_K}})$.
The next Theorem
follows from the results of Section~\ref{sect5} and from 
the multiplication formula of \cite{GLS3}.

\begin{Thm}
\begin{itemize}
\item[(i)]
$\{\x(T)\mid T\in \R_J\}$ is the set of clusters
of a cluster algebra ${\cal A}_J\subseteq \C[N_K]$
of rank $d_K-n$.
\item[(ii)]
The coefficient ring of ${\cal A}_J$ is the
ring of polynomials in the $n$ variables $\varphi_{L_i}\ (i\in I)$.
\item[(iii)]
All the cluster monomials belong to the dual semicanonical
basis of $\C[N_K]$, and are thus linearly independent.
\end{itemize}
\cqfd   
\end{Thm}

\subsection{}
We now proceed to describe in detail initial seeds for
${\cal A}_J$.

\subsubsection{}
Let $\ii=(i_1,\ldots,i_r)\in R(w_0,K)$.
Put $i_{-m}=m$ for $m=1,\ldots ,n$.
As in \ref{ss25}, set $u_{\le m}=s_{i_1}\cdots s_{i_m}$
if $m=1,\ldots,r$, and $u_{\le m}= e$ the unit element of $W$
if $m<0$.
Then, by Corollary~\ref{Cor62}, the cluster of $\C[N]$
obtained by projecting $T^\dag_\ii$ consists of the r functions
\begin{equation}
\varphi(m,\ii)=D_{u_{\le m}(\varpi_{i_m}),\,w_0(\varpi_{i_m})},
\qquad (m\in [-n,-1]\cup e(\ii)).
\end{equation}
Moreover, $T^\dag_\ii=ST_\ii$, where, as seen in  \ref{ss25}, 
$T_\ii$ projects on
one of the initial seeds of \cite{BFZ}. 
By \cite{GLS2} the exchange matrix of a cluster seed of $\C[N]$
is a submatrix of the matrix of the Ringel form of the endomorphism
algebra of the corresponding complete rigid module. 
It follows from \ref{ss24} that the exchange matrix of the seeds corresponding
to $T_\ii$ and $T_\ii^\dag$ are given by the same combinatorial
rule, which is described in \cite{BFZ}. 
We shall recall it for the convenience of the reader.

\subsubsection{}
For $m \in [-n,-1] \cup e(\ii)$, let 
$m^+ =  \min\{l\in[1,r] \mid l > m \mbox{ and } i_l = i_m \}$.
Next, one defines a quiver $\Gamma_\ii$ with set of vertices
$[-n,-1] \cup e(\ii)$.
Assume that $m$ and $l$ are vertices such that $m < l$ and
$\{ m,l \} \cap e(\ii) \not= \emptyset$.
There is an arrow $m \to l$ in $\Gamma_\ii$ if and only if
$m^+ = l$, 
and there is an arrow $l \to m$ if and only if
$l < m^+ < l^+$ and $a_{i_m,i_l} <0$.
Here, $(a_{ij})_{1\le i,j \le n}$ denotes the Cartan matrix
of the root system of $G$.
By definition these are all the arrows of $\Gamma_\ii$.
Now define an $r \times (r-n)$-matrix 
$B(\ii) = (b_{ml})$
as follows.
The columns of $B(\ii)$ are indexed by $e(\ii)$, 
and the rows by $e(\ii) \cup [-n,-1]$.
Set
\[
b_{ml} = \left\{
\begin{array}{ll}
1 & \mbox{if there is an arrow $m\to l$ in $\Gamma_\ii$,}\\[2mm]
-1 & \mbox{if there is an arrow $l \to m$ in  $\Gamma_\ii$,}\\[2mm]
0 &\mbox{otherwise.}
\end{array}
\right.
\]

\subsubsection{}\label{ss84}
Finally, as in \ref{ss52}, for $k\in K$ let $t_k=\max\{t\le r_K \mid i_t = k\}$.
If $j\in J$, set $t_j=-j$.
It follows from \ref{ss52} that 
\[
\x(U_\ii)=\{\varphi(m,\ii)\mid m\in ]r_K,r]\cap e(\ii)\}
\cup \{\varphi(t_i,\ii)\mid i\in I\}.
\]
This is the cluster of our initial seed for ${\cal A}_J$.
Here the first subset consists of the $r-r_K-n$ cluster
variables. The second subset consisting of the $n$ generators
of the coefficient ring is contained in every cluster. 
The corresponding exchange matrix is the submatrix of $B(\ii)$
with columns labelled by $]r_K,r]\cap e(\ii)$ and
rows labelled by $\left(]r_K,r]\cap e(\ii)\right)\cup \{t_i\mid i\in I\}$. 
Let us denote it by $B(\ii,J)$.
Then our construction yields
\begin{Prop}\label{prop94}
For every $\ii\in R(w_0,K)$, the pair $(\x(U_\ii),B(\ii,J))$
is an initial seed of the cluster algebra 
${\cal A}_J\subseteq \C[N_K]$.
\cqfd
\end{Prop}

\begin{example}\label{ex95}
{\rm
Let $G$ be of type $A_5$, that is, $G=SL_6$. Then $r=15$.
Take $J=\{1,3\}$. Then we have $K=\{2,4,5\}$ and $r_K=4$.
The word
$\ii=({\bf 2,4,5,4,} \ 1,2,3,4,5,2,3,4,1,2,3)$
belongs to $R(w_0,K)$ :
the subword $(2,4,5,4)$ is a reduced word for $w_0^K$.
We have 
\[
e(\ii)=\{1,2,3,4,5,6,7,8,10,11\}.
\]
\begin{figure}[t]
\begin{center}
\leavevmode
\epsfxsize =11cm
\epsffile{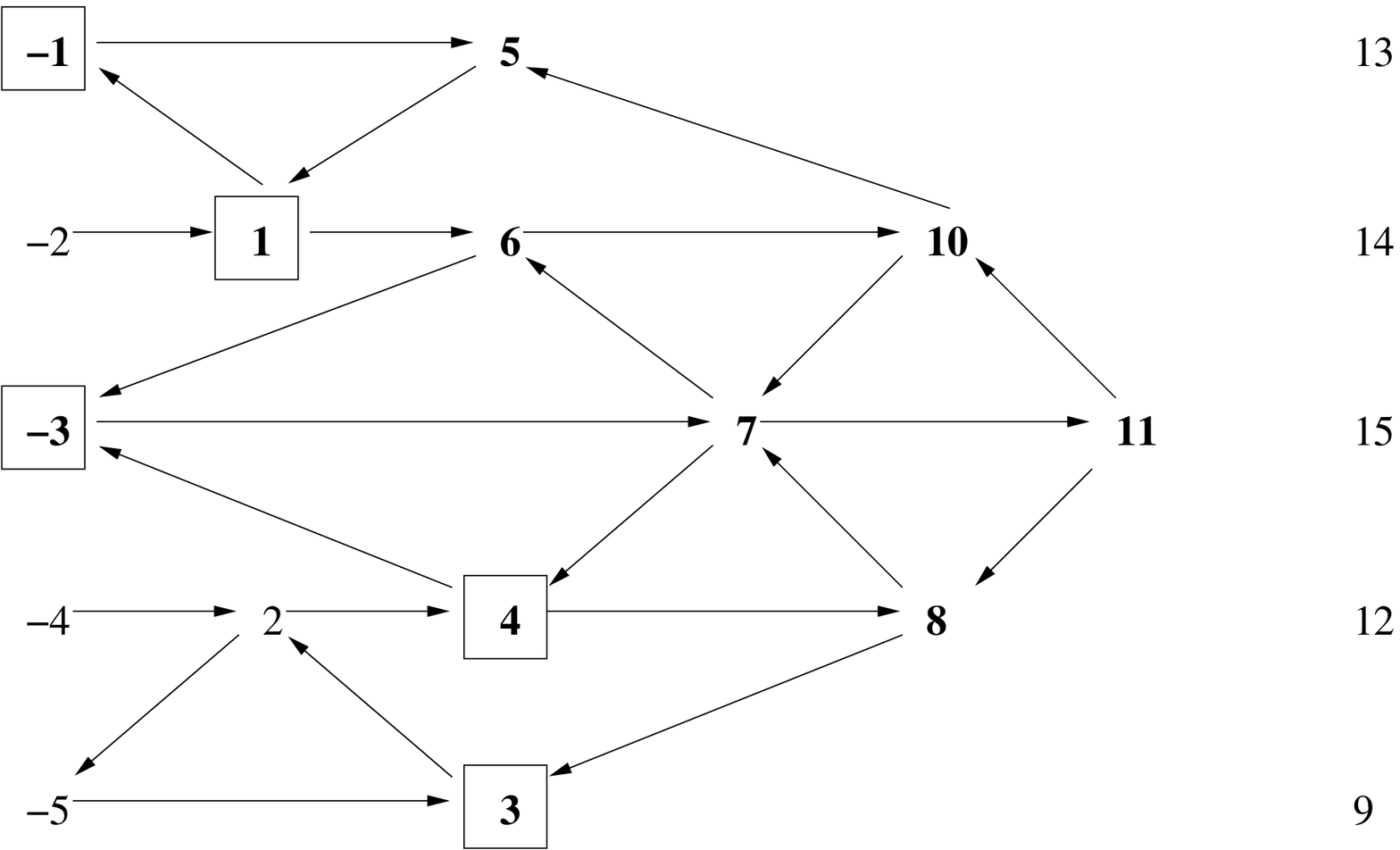}
\end{center}
\caption{\small {\it The graph $\Gamma_\ii$.}}\label{fig1}
\end{figure}
The graph $\Gamma_\ii$ is shown in Figure~\ref{fig1}.
The vertices $5, 6, 7, 8, 10, 11$ (in bold type) represent
the cluster variables of the initial seed of ${\cal A}_J$,
and the vertices $-1, 1, -3, 4, 3$ (in squares) 
represent the generators of the coefficient ring. 
The exchange matrix is 
\[
B(\ii,J)=
\begin{pmatrix}
         0&0&0&0&-1&0\cr
         0&0&-1&0&1&0\cr
         0&1&0&-1&-1&1\cr
         0&0&1&0&0&-1\cr
         1&-1&1&0&0&-1\cr
         0&0&-1&1&1&0\cr
         \hline
         1&0&0&0&0&0\cr
         -1&1&0&0&0&0\cr
         0&-1&1&0&0&0\cr
         0&0&-1&1&0&0\cr
         0&0&0&-1&0&0           
\end{pmatrix},
\]
where 
the rows of $B(\ii,J)$ are labelled by 
$\{5, 6, 7, 8, 10, 11,-1, 1, -3, 4, 3\}$ and the columns by 
$\{5, 6, 7, 8, 10, 11\}$.
We have
\[
\varphi(5,\ii) = D_{3,6},\quad
\varphi(6,\ii) = D_{23,56},\quad
\varphi(7,\ii) = D_{236,456},\quad
\varphi(8,\ii) = D_{2356,3456},\quad
\]
\[
\varphi(10,\ii) = D_{36,56},\quad
\varphi(11,\ii) = D_{356,456},\quad
\varphi(-1,\ii) = D_{1,6},\quad
\varphi(1,\ii) = D_{13,56},\quad
\]
\[
\varphi(-3,\ii) = D_{123,456},\quad
\varphi(4,\ii) = D_{1236,3456},\quad
\varphi(3,\ii) = D_{12356,23456}.
\]
Here, since we are in type $A_5$, the generalized
minors are simply minors of the upper unitriangular $6\times 6$ matrix of 
coordinate functions on $N$. 
More precisely, for $w$ in the symmetric group
$\SG_6$ the weight $w(\varpi_i)$ indexing
a generalized minor can be identified with the
subset $w([1,i])$ of $[1,6]$ of row or column indices
of the minor.
This is the convention we have used.
\finex 
}
\end{example}

\subsection{}
We end this section by stating the following

\begin{Conj}\label{conj1}
We have ${\cal A}_J = \C[N_K]$.
\end{Conj}
The conjecture will be proved for $G$ of type $A_n$ and of type $D_4$
in \ref{subsecconj2}.
It will also be proved for
$J=\{n\}$ in type $D_n$ in Section~\ref{sect12}, and for $J=\{1\}$
in type $D_5$ in Section~\ref{sect13}.


\section{Cluster algebra structure on $\C[B_K^-\backslash G]$}


In this section we lift the cluster algebra ${\cal A}_J\subseteq \C[N_K]$
to a cluster algebra $\widetilde{\cal A}_J \subseteq \C[B_K^-\backslash G]$.

\subsection{}
Let $(\x(T), B)$ denote a seed of ${\cal A}_J$.
Let $x_k=\varphi_{T_k}$ be a cluster variable in $\x(T)$, and denote by
$x_k^*$ the cluster variable obtained by mutation in
direction $k$. The exchange relation between $x_k$ and 
$x_k^*$ is of the form
\[
x_kx_k^*= M_k + N_k,
\]
where $M_k$ and $N_k$ are monomials in the variables of $\x(T)-\{x_k\}$.
We have two exact sequences in $\sub Q_J$
\[
0\to T_k \to X_k \to T_k^*\to 0,
\qquad 
0\to T_k^* \to Y_k \to T_k\to 0,
\]
where $X_k, Y_k\in \add T$,  and 
\[
x_k^*=\varphi_{T_k^*},\qquad M_k=\varphi_{X_k},\qquad N_k=\varphi_{Y_k}.
\] 
The following proposition is a particular case of
a result in \cite{GLS5}.
\begin{Prop}\label{prop10.1}
With the above notation,
for every $j\in J$ we have
\[
\dim \Hom_\L(S_j,T_k) + \dim \Hom_\L(S_j,T_k^*)
=
\max\{\dim \Hom_\L(S_j,X_k),\,\dim \Hom_\L(S_j,Y_k)\}.
\] \cqfd
\end{Prop}
Recall from Section~\ref{ss93} the definition of $a_j(f)$
for $f\in \C[N_K]$.
If $f=\varphi_M$ for some $\L$-module~$M$, it follows
from \cite{GLS1} that $a_j(f)=\dim \Hom_\L(S_j,M)$,
the multiplicity of $S_j$ in the socle of $M$.
Putting together Proposition~\ref{prop10.1} and 
Lemma~\ref{lem25} we thus have
\begin{equation}\label{eq94}
\widetilde{x_k}\widetilde{x_k^*}= 
\mu_k\widetilde{M_k} + \nu_k\widetilde{N_k},
\end{equation}
where $\mu_k$ and $\nu_k$ are monomials in the variables
$\De_{\varpi_j,\varpi_j}\ (j\in J)$ without common divisor.
More precisely,
\begin{equation}\label{eq8}
\mu_k = \prod_{j\in J} \De_{\varpi_j,\varpi_j}^{\alpha_j},
\qquad
\nu_k = \prod_{j\in J} \De_{\varpi_j,\varpi_j}^{\beta_j},  
\end{equation}
where
\begin{equation}\label{eq9}
\alpha_j = \max\{0\,,\,\dim\Hom_\L(S_j,Y_k)-\dim \Hom_\L(S_j,X_k)\},
\end{equation}
\begin{equation}\label{eq10}
\beta_j = \max\{0\,,\,\dim\Hom_\L(S_j,X_k)-\dim \Hom_\L(S_j,Y_k)\}.
\end{equation}
\subsection{}\label{ssect10.2}
It follows that the elements $\tx$ where $x$ runs through
the set of cluster variables of ${\cal A}_J$ form the
cluster variables of a new cluster algebra 
contained in $\C[B_K^-\backslash G]$.
More precisely, for
$T=T_1\oplus \cdots \oplus T_{d_K}\in\R_J$, 
let $\widetilde{\x(T)}=
\{\widetilde{\varphi_{T_1}},\ldots,\widetilde{\varphi_{T_{d_K}}}\}
\sqcup \{\De_{\varpi_j,\varpi_j}\mid j\in J\}$.

\begin{Thm}
\begin{itemize}
\item[(i)]
$\{\widetilde{\x(T)}\mid T\in \R_J\}$ is the set of clusters
of a cluster algebra $\widetilde{\cal A}_J\subseteq \C[B_K^-\backslash G]$
of rank $d_K-n$.
\item[(ii)]
The coefficient ring of $\widetilde{\cal A}_J$ is the
ring of polynomials in the $n+|J|$ variables 
$\widetilde{\varphi_{L_i}}\ (i\in I)$
and $\De_{\varpi_j,\varpi_j}\ (j\in J)$.
\item[(iii)]
The exchange matrix $\tB$ attached to $\widetilde{\x(T)}$
is obtained from the exchange matrix $B$ of $\x(T)$ by
adding $|J|$ new rows (in the non-principal part)
labelled by $j\in J$
encoding the monomials $\mu_k$ and $\nu_k$ in Equations~(\ref{eq94})
and (\ref{eq8}).
More precisely, the entry in column $k$ and row $j$ is
equal to 
\[
b_{jk}=\dim\Hom_\L(S_j,X_k)-\dim \Hom_\L(S_j,Y_k).
\]
\end{itemize}
\end{Thm}
\proof
The only thing to be proved is that the enlarged exchange
matrices $\tB$ of (iii) follow the matrix mutation rule of
Fomin and Zelevinsky.
To do that, let us introduce some notation.
Let $b_{ik}$ denote the entry of $\tB=\tB(T)$ on row $i$ and
column $k$. Here $k$ runs through $\{1,\ldots, d_K-n\}$
and $i$ runs through $\{1,\ldots, d_K\}\cup J$.
For $k\in\{1,\ldots, d_K-n\}$ the mutation of $T$ in direction $k$ 
comes from the two short exact sequences
\[
0\to T_k \to X_k \to T_k^*\to 0,
\qquad 
0\to T_k^* \to Y_k \to T_k\to 0.
\]
Fix $l\in\{1,\ldots, d_K-n\}$ and let us consider the mutation
$\mu_l(T)=T_l^*\oplus T/T_l$.
We denote by $\tB'=[b_{ik}']$ the exchange matrix attached
to $\mu_l(T)$, as defined in (iii).
We have to show that $b_{ik}'=-b_{ik}$ if $i$ or $k$ is equal
to $l$, and otherwise 
\[
b'_{ik}=b_{ik}+{b_{il}|b_{lk}|+|b_{il}|b_{lk}\over 2}.
\]
(a) By construction this holds for $i\in \{1,\ldots, d_K\}$, so
we shall now assume that $i$ is an element $j$ of~$J$.

If $k=l$, the two short exact sequences used to define $b_{jk}$
and $b'_{jk}$ are the same but they are interchanged, so it
is plain that $b'_{jk}=-b_{jk}$. 

From now on we assume that $k\not = l$.
Let
\[
0\to T_k \to X_k' \to {}'T_k^*\to 0,
\qquad 
0\to {}'T_k^* \to Y_k' \to T_k\to 0,
\]
denote the two short exact sequences inducing the mutation of
$\mu_l(T)$ in direction $k$.
If $b_{lk}=0$ then $b'_{ik}=b_{ik}$ for all
$i\in\{1,\ldots, d_K\}$. It follows that $X'_k=X_k$ and
$Y'_k=Y_k$, hence ${}'T_k^*=T_k^*$, $b'_{jk}=b_{jk}$,
and the result is proved. 
Hence we can assume that $b_{lk}\not = 0$.

Assume that $b_{lk} > 0$.
For $X,Y\in\mod\L$ we shall use the shorthand notation
\[
[X,Y]=\dim\Hom_\L(X,Y).
\]
We have 
\[
[S_j,X_k]-[S_j,Y_k]=-\sum_{b_{ik}<0} b_{ik} [S_j,T_i]
-\sum_{b_{ik}>0} b_{ik} [S_j,T_i] = -\sum_i b_{ik} [S_j,T_i],
\]
where the sums are over $i\in\{1,\ldots, d_K\}$.
Now, for $i\in\{1,\ldots, d_K\}$, it follows from (a) that
$b'_{ik}= b_{ik}$ unless $b_{il}>0$ in which case
$b'_{ik}= b_{ik}+b_{il}b_{lk}$.

(b) Let us further assume that $b_{jl}>0$. This implies
that $[S_j,T_l^*]=[S_j,X_l]-[S_j,T_l]$. Now, using
that
\[
X'_k=(T_l^*)^{b_{lk}}\oplus
\bigoplus_{\substack{i\not = l\\b'_{ik}<0}} T_i^{(-b'_{ik})},
\qquad
Y'_k=\bigoplus_{b'_{ik}>0} T_i^{b'_{ik}},
\]
we calculate
\[
[S_j,X'_k]=b_{lk}\left([S_j,X_l]-[S_j,T_l]\right)
-\sum_{\substack{b_{il}>0\\b_{ik}+b_{il}b_{lk}<0}} 
(b_{ik}+b_{il}b_{lk})[S_j,T_i]
-\sum_{\substack{b_{il}\le0\\b_{ik}<0}} b_{ik}[S_j,T_i],
\]
\[
[S_j,Y'_k]=
\sum_{\substack{b_{il}>0\\b_{ik}+b_{il}b_{lk}>0}} 
(b_{ik}+b_{il}b_{lk})[S_j,T_i]
+\sum_{\substack{b_{il}\le0\\b_{ik}>0\\i\not=l}} b_{ik}[S_j,T_i].
\]
This implies that
\[
[S_j,X'_k]-[S_j,Y'_k]=
b_{lk}([S_j,X_l]-[S_j,T_l])
-\sum_{b_{il}>0} 
(b_{ik}+b_{il}b_{lk})[S_j,T_i] 
-\sum_{\substack{b_{il}\le0\\i\not=l}} b_{ik}[S_j,T_i]
\]
\[
=b_{lk}[S_j,X_l]
-\sum_{b_{il}>0} 
(b_{ik}+b_{il}b_{lk})[S_j,T_i] 
-\sum_{b_{il}\le0} b_{ik}[S_j,T_i]
\]
\[
=
-\sum_{i} b_{ik}[S_j,T_i]
+b_{lk}\left([S_j,X_l]-\sum_{b_{il}>0}b_{il}[S_j,T_i]\right).
\]
Hence
\[
[S_j,X'_k]-[S_j,Y'_k]=[S_j,X_k]-[S_j,Y_k]+b_{lk}([S_j,X_l]-[S_j,Y_l]),
\]
that is $b'_{jk}=b_{jk}+b_{lk}b_{jl}$, as required.

(c) Assume now that $b_{jl} \le 0$.
This implies that $[S_j,T_l^*]=[S_j,Y_l]-[S_j,T_l]$,
and the same calculation as above now gives
\[
[S_j,X'_k]-[S_j,Y'_k]=[S_j,X_k]-[S_j,Y_k],
\]
that is, $b'_{jk}=b_{jk}$.
This finishes the proof when $b_{lk} > 0$. The case $b_{lk} < 0$
is entirely similar.
\cqfd

Note that $d_K+|J|$ is equal to the dimension of the
multi-cone over $B_K^-\backslash G$.
Note also that the clusters of ${\cal A}_J$ and  $\widetilde{\cal A}_J$
are in natural one-to-one correspondence, and the principal
parts of the exchange matrices of two corresponding clusters
are the same. This shows that ${\cal A}_J$ and $\widetilde{\cal A}_J$
are both of finite cluster type or of infinite cluster type,
and if they are of finite type, their types are the same.
\begin{example}\label{ex10.2}
{\rm
We continue Example~\ref{ex95}.
Let us denote for short by $x_m$ the functions
$\varphi(m,\ii)$ of Example~\ref{ex95}.
Thus we have $11$ minors 
$$
x_5,x_6,x_7, x_8, x_{10}, x_{11},x_{-1}, x_1, x_{-3}, x_4,x_3,
$$
in $\C[N_K]$.
It is straightforward to calculate their lifts to $\C[B_K^-\backslash G]$. 
One gets
 \[
\widetilde{x_5} = \De_{126},\quad
\widetilde{x_6} = \De_{156},\quad
\widetilde{x_7} = \De_{145},\quad
\widetilde{x_8} = \De_{134},\quad
\widetilde{x_{10}} = \De_{125},\quad
\widetilde{x_{11}} = \De_{124},
\]
\[
\widetilde{x_{-1}} = \De_6,\quad
\widetilde{x_1} = \De_2\De_{156}-\De_1\De_{256},\quad
\widetilde{x_{-3}} = \De_{456},\quad
\widetilde{x_4} = \De_{345},\quad
\widetilde{x_3} = \De_{234}.
\]
In these formulas, all minors are flag minors, hence
we have indicated only their column indices.
For example, using the fact that $D_{I,J}$ is a minor of an upper
unitriangular matrix, we have, 
\[
x_8 = D_{2356,3456}
    = D_{23,34}
    = D_{123,134}.
\]
Hence, $\widetilde{x_8} = \De_{123,134}= \De_{134}$.
A more interesting example is 
\[
x_1 = D_{13,56}
    = D_{1, 2}D_{23,56}-D_{123,256}
    = D_{1, 2}D_{123,156}-D_{123,256}.
\]
Hence, 
\[
\widetilde{x_1} = \De_{1, 2}\De_{123,156}-\De_{1,1}\De_{123,256} 
=\De_2\De_{156}-\De_1\De_{256}.
\]
Note that $\widetilde{x_1}$ cannot be written as a flag minor on $G$.
(In representation theoretical terms, this is because the
socle of the indecomposable rigid $\L$-module $L_4$ attached to
$x_1$ is not simple : it is isomorphic to $S_1\oplus S_3$.)
Since $\widetilde{x_1}$ is one of the generators of the coefficient 
ring of $\widetilde{\cal A}_J$
(that is, $L_4=P_4/S_2$ is projective in $\sub Q_J$, see Proposition~\ref{prop2-1}), 
it belongs to every seed of the
cluster structure of $\widetilde{\cal A}_J$.
This example shows that, in contrast to the case of the full flag
variety, the cluster algebra $\widetilde{\cal A}_J$ of a partial
flag variety may have no seed consisting entirely of flag minors. 

Finally, the exchange matrix for this seed of 
$\widetilde{\cal A}_J$ is
\[
\tB(\ii,J)=
\begin{pmatrix}
         0&0&0&0&-1&0\cr
         0&0&-1&0&1&0\cr
         0&1&0&-1&-1&1\cr
         0&0&1&0&0&-1\cr
         1&-1&1&0&0&-1\cr
         0&0&-1&1&1&0\cr
         \hline        
         1&0&0&0&0&0\cr
         -1&1&0&0&0&0\cr
         0&-1&1&0&0&0\cr
         0&0&-1&1&0&0\cr
         0&0&0&-1&0&0\cr
         \hline
         0&-1&0&0&0&0\cr
         0&0&0&0&0&1
\end{pmatrix}.
\]
It is obtained from the matrix of Example~\ref{ex95}
by adding the last two rows labelled by $\De_1$ and $\De_{123}$.
The two nonzero entries in these rows correspond to the
two exchange relations for $\widetilde{x_{6}}$ 
and~$\widetilde{x_{11}}$.
\finex
}
\end{example}

\subsection{}
We now show that in type $A_n$ our algebras $\widetilde{\cal A}_{\{j\}}$
coincide with the cluster algebras on coordinate rings of Grassmannians
considered in \cite{FZ3} for $j=2$, and in \cite{S} for general $j$.
To do so, it is enough to check that one of our seed for 
$\widetilde{\cal A}_{\{j\}}$ coincides with a seed of \cite{S}.

\subsubsection{}\label{grid}
In the case of type $A_n$ and $J=\{j\}$, 
the rule described in Proposition~\ref{prop94} gives us a
unique initial cluster for ${\cal A}_{\{j\}}$. Indeed, if 
$\ii=(\ii',\ii'')$ and $\jj=(\jj',\jj'')$ are 
two elements in $R(w_0,I-\{j\})$, then $\ii''$ and
$\jj''$ are related by a sequence of $2$-moves, and
the corresponding clusters are therefore equal.
This unique initial cluster of ${\cal A}_{\{j\}}$ consists 
of the minors $D_C$ where $C$ belongs to the following list:
\[
\{1,\ldots,j-1,j+1\},\
\{1,\ldots,j-1,j+2\},\
\ldots, 
{\bf\{1,\ldots,j-1,n+1\}},
\]
\[
\{1,\ldots,j-2,j,j+1\},\
\{1,\ldots,j-2,j+1,j+2\},\
\ldots, 
{\bf\{1,\ldots,j-2,n,n+1\}},
\]
\[
\{1,\ldots,j{-}3,j{-}1,j,j{+}1\},\
\{1,\ldots,j{-}3,j,j{+}1,j{+}2\},\
\ldots, 
{\bf\{1,\ldots,j{-}3,n{-}1,n,n{+}1\}},
\]
\[
\ldots,\qquad\qquad\ldots,\qquad\qquad\ldots,
\]
\[
\{1,3,\ldots,j,j+1\},\
\{1,4,\ldots,j+1,j+2\},\
\ldots, 
{\bf\{1,n-j+3\ldots,n,n+1\}},
\]
\[
{\bf\{2,\ldots,j+1\},\
\{3,\ldots,j+2\},\
\ldots, 
\{n-j+2,\ldots,n+1\}}.
\]
The $n$ subsets in bold type correspond to generators of the
coefficient ring, and therefore cannot be mutated.
Here we use the same notation as in Example~\ref{ex10.2}
for the flag minors in type $A$.
Note that it is not obvious from Proposition~\ref{prop94} that
all elements of this initial cluster can be written as flag 
minors $D_C$. It turns out to be true and not difficult to
check in this special case. 
From the point of view of preprojective algebras, this case
is special in the sense that all indecomposable summands of
the rigid module $U_\ii$ of \ref{sect9.2}, which by construction 
have simple top, also have simple socle.

The graph describing the exchange matrix of this cluster
has the shape of a rectangular grid with $j$ rows and
$n-j+1$ columns. The vertices are the subsets $C$ 
displayed as above. There are horizontal right arrows,
vertical down arrows and diagonal north-west arrows.
For example, if $n=7$ and $j=4$, this graph is shown  
\begin{figure}[t]
\begin{center}
\leavevmode
\epsfxsize =8cm
\epsffile{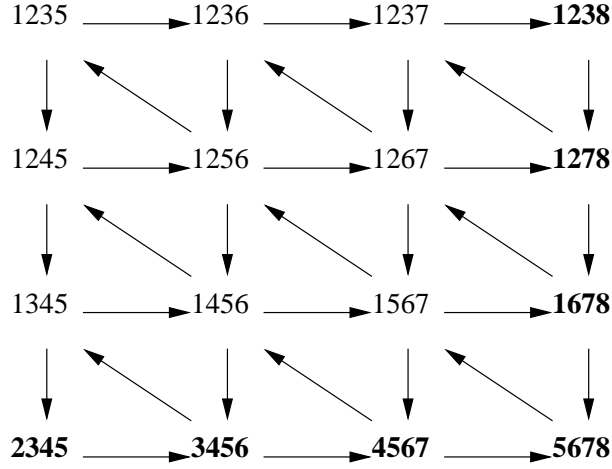}
\end{center}
\caption{\small {\it The grid for $n=7$ and $j=4$.}}\label{fig2}
\end{figure}
in Figure~\ref{fig2}.

We note that this initial seed coincides with the one described in 
\cite[\S 3.3]{GSV}.

\subsubsection{}
\begin{figure}[t]
\begin{center}
\leavevmode
\epsfxsize =12cm
\epsffile{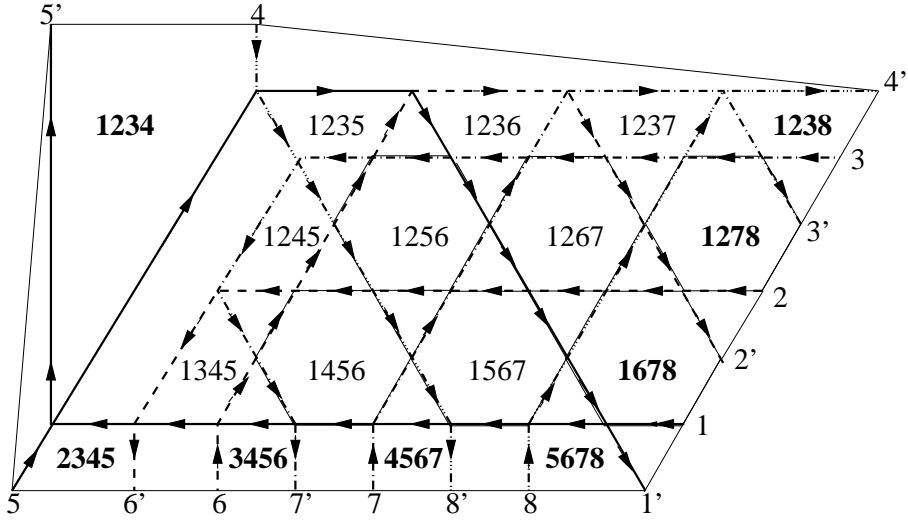}
\end{center}
\caption{\small {\it The honeycomb arrangement for $n=7$ and $j=4$.}}
\label{honey}
\end{figure}

Let $\De_C$ denote the flag minor on $G$ corresponding
to $D_C$. Thus, $\De_C=\widetilde{D_C}$.
The initial cluster of $\widetilde{\cal A}_{\{j\}}$
lifting the cluster above consists of all $\De_C$
where $C$ runs over all sets in the above list, together
with the new set $[1,j]=\{1,\ldots,j\}$.
All the exchange relations are obtained by simply replacing
each $D_C$ by the corresponding $\De_C$, except the exchange
relation for $C=\{1,\ldots,j-1,j+1\}$ which reads
\[
\De_{\{1,\ldots,j{-}1,j{+}1\}} \De_{\{1,\ldots,j{-}2,j,j{+}2\}} 
=
\De_{\{1,\ldots,j{-}1,j{+}2\}} \De_{\{1,\ldots,j{-}2,j,j{+}1\}} 
+
\De_{[1,j]} \De_{\{1,\ldots,j{-}2,j{+}1,j{+}2\}}.
\]
One then checks that this coincides with one of the 
seeds of \cite{S}.
To do so, one has to construct a Postnikov arrangement \cite{P},
\cite[\S 3]{S}
whose labelling is given by the list of subsets $C$ of our seed
for $\widetilde{\cal A}_{\{j\}}$.
This arrangement has a regular structure similar to a honeycomb.
The central labelled cells are hexagonal and the border ones
are quadrilateral.
If the sets $C$ and $C'$ are connected by an arrow in the grid
of \ref{grid}, then the cells labelled by $C$ and $C'$ have a common
vertex. 
For example, when $n=7$ and $j=4$, the ``honeycomb'' Postnikov arrangement
is shown in Figure~\ref{honey}.
This finishes the proof that $\widetilde{\cal A}_{\{j\}}$ coincides
with the cluster algebra of \cite{S}.

\subsubsection{}
\begin{figure}[t]
\begin{center}
\leavevmode
\epsfxsize =14.5cm
\epsffile{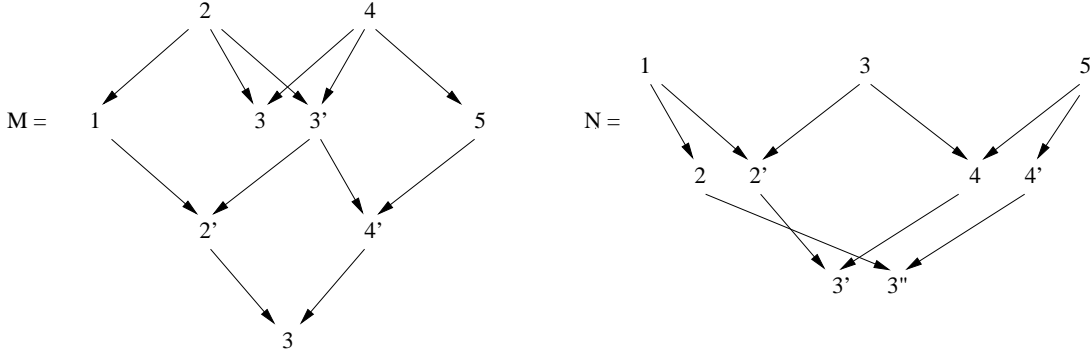}
\end{center}
\caption{\small {\it The rigid $\L$-modules $M$ and $N$ in $\sub Q_3$ 
for type $A_5$.}}
\label{fig3}
\end{figure}

We think that our construction sheds a new light on the cluster
algebra structures of the coordinate rings of Grassmannians.
For example, in \cite[Th. 6]{S}, Scott describes the two 
special cluster variables $X$ and $Y$ of $\C[\Gr(3,6)]$ which are not 
flag minors. In our setting 
\[
X=\widetilde{\varphi_M}, \qquad Y=\widetilde{\varphi_N},
\]
where $M$ and $N$ are the only two indecomposable rigid
modules of $\sub Q_3$ in type $A_5$ with a $2$-di\-men\-sio\-nal 
socle $S_3\oplus S_3$.
These modules are represented in Figure~\ref{fig3}.

\subsection{}\label{subsecconj2}
Let $\Si_J$ be the multiplicative submonoid of $\widetilde{\cal A}_J$
generated by the set
\[
\{\De_{\varpi_j,\varpi_j} \mid j\in J \mbox{ and } \varpi_j
\mbox{ is {\em not } minuscule}\}.
\]

\begin{Conj}\label{conj2}
The localizations of $\widetilde{\cal A}_J$ and 
$\C[B_K^-\backslash G]$ with respect to $\Si_J$ are equal.
\end{Conj}
Note that if $J$ is such that all the weights
$\varpi_j\ (j\in J)$ are minuscule, then $\Si_J$ is trivial
and the conjecture states that the algebras $\widetilde{\cal A}_J$ and 
$\C[B_K^-\backslash G]$ coincide without localization.
This is in particular the case for every $J$ in type $A_n$.
Note also that Conjecture~\ref{conj2} implies 
Conjecture~\ref{conj1}.
Indeed, by construction 
$\pr_J(\widetilde{\cal A}_J)= {\cal A}_J$,
$\pr_J(\Si_J) = \{1\}$
and by \ref{ss93},  
$\pr_J(\C[B_K^-\backslash G])=\C[N_K]$.
We shall now prove Conjecture~\ref{conj2} in type $A_n$ and $D_4$.

\subsubsection{}\label{conjgrass}
We first remark that if the conjecture is true for 
every $\widetilde{\cal A}_{\{j\}}\ (j\in J)$, then it is true for
$\widetilde{\cal A}_J$.
Indeed, if $J'\subset J$ and $K'=I\setminus J$,
there exist reduced words for $w_0$ of the form
$\ii=(\ii',\ii'',\ii''')$ with 
$\ii'$ a reduced word for $w_0^{K}$ and $(\ii',\ii'')$
a reduced word for $w_0^{K'}$. This shows that the 
initial seed for $\widetilde{\cal A}_J$ associated with $\ii$ will
contain the initial seed for $\widetilde{\cal A}_{J'}$ 
associated with $\ii$. Hence $\widetilde{\cal A}_{J'}$ is a 
subalgebra of $\widetilde{\cal A}_J$. In particular, 
$\widetilde{\cal A}_{\{j\}}$ is a subalgebra of $\widetilde{\cal A}_J$ 
for every $j\in J$.
Suppose we know that 
$\C[B^-_{I\setminus\{j\}}\backslash G]$
is contained in the localization of
$\widetilde{\cal A}_{\{j\}}$ with respect to $\Si_{\{j\}}$
for every $j\in J$.
Then, since $\C[B^-_K\backslash G]=\bigoplus_{\l\in\Pi_J}L(\l)$ 
is generated (as a ring) by the subspaces 
$L(\varpi_j)\subset \C[B^-_{I\setminus\{j\}}\backslash G]$ 
we have that $\C[B^-_K\backslash G]$ is generated by the 
localized rings $\widetilde{\cal A}_{\{j\}}[\Si_{\{j\}}^{-1}]$,
hence also by $\widetilde{\cal A}_J[\Si_J^{-1}]$. 
Therefore Conjecture~\ref{conj2} is satisfied.  

\subsubsection{}\label{sss942}
In type $A_n$, the algebra $\C[B^-_{I\setminus\{j\}}\backslash G]$ 
is generated by the set of flag minors 
\[
\De_{\varpi_j,w(\varpi_j)}\qquad (w\in W).
\]
Hence to prove the conjecture in this case it is enough
to show that each of these minors belongs to $\widetilde{\cal A}_{\{j\}}$.
In \cite{S}, Scott proves that all flag minors are cluster variables
of $\widetilde{\cal A}_{\{j\}}$. It follows that 
Conjecture~\ref{conj2} is true in type $A$.

\subsubsection{}
We now prove Conjecture~\ref{conj2} in type $D_4$.
We choose to label by 3 the central node of the Dynkin diagram.
By \ref{conjgrass} it is enough to check the conjecture in the cases 
$J=\{1\}, \{2\}, \{3\}$ and $\{4\}$.
Because of the order 3 diagram automorphism of $D_4$,
the cases $J=\{1\}, \{2\}$ and $\{4\}$ are identical.
They are dealt with in detail in Section~\ref{sect12}
(which studies more generally the case $J=\{n\}$ in type
$D_n$). 
So we are left with $J=\{3\}$.

\begin{figure}[t]
\begin{center}
\leavevmode
\epsfxsize =3.5cm
\epsffile{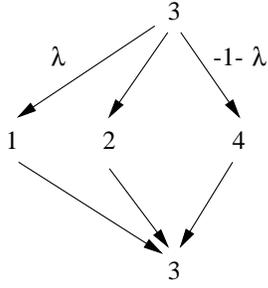}
\end{center}
\caption{\small {\it The one-parameter family $M(\l)$ in $\sub Q_3$ 
for type $D_4$.}}
\label{figM}
\end{figure}

To prove the conjecture in this case we have to show
that $\widetilde{\cal A}_{\{3\}}$ contains, up to localization
by $\De_{\varpi_3,\varpi_3}$, a basis of 
the $G$-module $L(\varpi_3)$.
This module has dimension 28.
The 24 generalized minors 
$\De_{\varpi_3,u(\varpi_3)}\ (u\in W)$
are extremal vectors of $L(\varpi_3)$, but we need
4 more vectors to get a basis.
We shall use the dual semicanonical basis of $L(\varpi_3)$,
obtained by lifting the dual semicanonical basis
of the subspace $\pr_{\{3\}}(L(\varpi_3))$ of $\C[N]$ via the
map $x \mapsto \widetilde{x}$.
This basis consists of elements of the form
$\widetilde{\varphi_M}$, where $M$ runs through 
``generic'' submodules of $Q_3$.
It contains the 24 minors above, attached to 24
rigid submodules $M$. 
In particular, 
\[
\De_{\varpi_3,\varpi_3} = \widetilde{\varphi_{\mathbf{0}}},\qquad
\De_{\varpi_3,w_0(\varpi_3)} = \widetilde{\varphi_{Q_3}},
\]
where $\mathbf{0}$ denotes the trivial submodule of $Q_3$.
The 4 remaining vectors are 
$\widetilde{\varphi_{L_1}}$,
$\widetilde{\varphi_{L_2}}$,
$\widetilde{\varphi_{L_4}}$
and $\widetilde{\varphi_{M(\l)}}$, where 
\[
M(\l),\qquad (\l \in \C-\{0,1\}),
\]
denotes the $1$-parameter family of $\L$-modules represented in 
Figure~\ref{figM}.
Note that, by construction, $\widetilde{\varphi_{L_1}}$,
$\widetilde{\varphi_{L_2}}$,
$\widetilde{\varphi_{L_4}}$,
$\De_{\varpi_3,\varpi_3}$
and $\De_{\varpi_3,w_0(\varpi_3)} =\widetilde{\varphi_{Q_3}}$
belong to $\widetilde{\cal A}_{\{3\}}$, since they
are the generators of the coefficient ring.

We choose $\ii = (1,2,4,3,1,2,4,3,1,2,4,3) \in R(w_0,K)$ and
consider the corresponding initial seed of 
$\widetilde{\cal A}_{\{3\}}$.
The 5 cluster variables of this seed are
\[
z_1=\De_{\varpi_3,s_3(\varpi_3)},\quad
z_2=\De_{\varpi_3,s_1s_3(\varpi_3)},\quad
z_3=\De_{\varpi_3,s_2s_3(\varpi_3)},\quad
z_4=\De_{\varpi_3,s_4s_3(\varpi_3)},
\]
and a variable $z_5$ of degree $2\varpi_3$.
One can then obtain by mutation the 18 other minors
of the form $\De_{\varpi_3,u(\varpi_3)}$, as shown 
in Table~\ref{tablemin}.
\begin{table}
\begin{center}
\begin{tabular}
{|c|c|}
\hline
$u$ & mutation sequence for $\De_{\varpi_3,u(\varpi_3)}$\\
\hline
$s_2s_1s_3$ & $\mu_4\mu_1$\\
$s_4s_1s_3$ & $\mu_3\mu_1$\\ 
$s_4s_1s_3$ & $\mu_2\mu_1$\\
$s_3s_2s_1s_3$ & $\mu_4$\\
$s_3s_4s_1s_3$ & $\mu_3$\\
$s_3s_4s_2s_3$ & $\mu_2$\\
$s_1s_4s_2s_3$ & $\mu_4\mu_5\mu_1\mu_3\mu_2\mu_1$\\
$s_4s_3s_2s_1s_3$ & $\mu_5\mu_1\mu_3\mu_2\mu_1$\\
$s_2s_3s_4s_1s_3$ & $\mu_5\mu_1\mu_2\mu_4\mu_1$\\
$s_1s_3s_4s_2s_3$ & $\mu_5\mu_1\mu_4\mu_3\mu_1$\\
$s_3s_1s_4s_2s_3$ & $\mu_5\mu_4\mu_3\mu_2$\\
$s_4s_3s_1s_4s_2s_3$ & $\mu_4\mu_5\mu_4\mu_3\mu_2$\\
$s_2s_3s_1s_4s_2s_3$ & $\mu_3\mu_5\mu_4\mu_3\mu_2$\\
$s_1s_3s_1s_4s_2s_3$ & $\mu_2\mu_5\mu_4\mu_3\mu_2$\\
$s_2s_4s_3s_1s_4s_2s_3$ & $\mu_2\mu_5$\\
$s_1s_2s_3s_1s_4s_2s_3$ & $\mu_3\mu_5$\\
$s_1s_4s_3s_1s_4s_2s_3$ & $\mu_4\mu_5$\\
$s_1s_2s_4s_3s_1s_4s_2s_3$ & $\mu_1\mu_5$\\
\hline
\end{tabular}
\end{center}
\caption{\small \it Generalized minors $\De_{\varpi_3,u(\varpi_3)}$
obtained by cluster mutation.
\label{tablemin}}
\end{table} 
Here, for example, the first row means that $\De_{\varpi_3,s_2s_1s_3(\varpi_3)}$  
is the new cluster variable obtained by applying to the initial
cluster  a mutation $\mu_1$ with respect to the first
variable, followed by a mutation $\mu_4$ with respect
to the fourth variable. 

It remains to show that $\widetilde{\varphi_{M(\l)}}$ also belongs
to $\widetilde{\cal A}_{\{3\}}$, up to division by $\De_{\varpi_3,\varpi_3}$.
For this, we use the multiplication formula for 
the functions $\varphi_M$ of \cite{GLS3}. 
Let $N_1$ and $N_2$ denote the $\L$-modules
\[
N_1=
\def\objectstyle{\scriptstyle} \xymatrix@-1.2pc 
{1\ar[rd]&\ar[d]2&\ar[ld]4\\&3&},
\qquad
N_2=
\def\objectstyle{\scriptstyle} \xymatrix@-1.2pc 
{1\ar[rd]&&\ar[ld]2\ar[rd]&&\ar[ld]4\\&3&&3}.
\]
It is easy to check that 
$\widetilde{\varphi_{N_1}}=\De_{\varpi_3,s_1s_4s_2s_3(\varpi_3)}$
and that
$\widetilde{\varphi_{N_2}}$ is the new cluster 
variable obtained from the initial cluster of 
$\widetilde{\cal A}_{\{3\}}$ via one mutation with respect to $z_1$.
By \cite{GLS3} we have
\[
\varphi_{M(\l)} =  \varphi_{S_3}\varphi_{N_1}
-\varphi_{N_2}-\varphi_{L_1}-\varphi_{L_2}-\varphi_{L_4}.
\]
This implies that
\[
\De_{\varpi_3,\varpi_3}\widetilde{\varphi_{M(\l)}} =  
\widetilde{\varphi_{S_3}}\widetilde{\varphi_{N_1}}
-\widetilde{\varphi_{N_2}}-\De_{\varpi_3,\varpi_3}
(\widetilde{\varphi_{L_1}}
+\widetilde{\varphi_{L_2}}
+\widetilde{\varphi_{L_4}}),
\]
which shows that $\widetilde{\varphi_{M(\l)}}$ belongs
to the localization of $\widetilde{\cal A}_{\{3\}}$ with
respect to $\De_{\varpi_3,\varpi_3}$.
This finishes the proof of Conjecture~\ref{conj2}
in type $D_4$.


\section{Finite type classification}
\label{sect11}

\begin{table}
\begin{center}
\begin{tabular}
{|c|c|c|}
\hline
Type of $G$ & $J$ & Type of ${\cal A}_J$\\
\hline
$A_n$ $(n\ge 2)$& $\{1\}$ & --- \\
$A_n$ $(n\ge 2)$& $\{2\}$ & $A_{n-2}$ \\
$A_n$ $(n\ge 2)$& $\{1,2\}$ & $A_{n-1}$\\
$A_n$ $(n\ge 2)$& $\{1,n\}$ & $(A_1)^{n-1}$\\
$A_n$ $(n\ge 3)$& $\{1,n-1\}$ & $A_{2n-4}$ \\
$A_n$ $(n\ge 3)$& $\{1,2,n\}$ & $A_{2n-3}$ \\
\hline
$A_4$ & $\{2,3\}$ & $D_4$\\
$A_4$ & $\{1,2,3\}$ & $D_5$\\
$A_4$ & $\{1,2,3,4\}$ & $D_6$\\
\hline
$A_5$ & $\{3\}$ & $D_4$ \\
$A_5$ & $\{1,3\}$ & $E_6$\\
$A_5$ & $\{2,3\}$ & $E_6$\\
$A_5$ & $\{1,2,3\}$ & $E_7$\\
\hline
$A_6$ & $\{3\}$ & $E_6$ \\
$A_6$ & $\{2,3\}$ & $E_8$\\
\hline
$A_7$ & $\{3\}$ & $E_8$ \\
\hline
$D_n$ $(n\ge 4)$ & $\{n\}$ & $(A_1)^{n-2}$\\
\hline
$D_4$ & $\{1,2\}$ & $A_5$\\
\hline
$D_5$ & $\{1\}$ & $A_5$\\
\hline
\end{tabular}
\end{center}
\caption{\small \it Algebras ${\cal A}_J$ of finite cluster type.
\label{tab1}}
\end{table} 

\subsection{}
Using the explicit initial seed described in \ref{ss84}
it is possible to give a complete list of the
algebras ${\cal A}_J$ which are of finite type as 
cluster algebras.
The results are summarized in Table~\ref{tab1}.
Here, we label the vertices of the Dynkin diagram of type
$D_n$ in such a way that
$L(\varpi_1)$ and $L(\varpi_2)$ are the two spin representations 
and $L(\varpi_n)$ is the vector representation.
We have only listed one representative of each orbit under a
diagram automorphism. For example, in type $A_n$
we have an order $2$ diagram automorphism mapping
$J=\{1,2\}$ to $J'=\{n-1,n\}$. Clearly,
${\cal A}_{J'}$ has the same cluster type as $J$,
namely $A_{n-1}$. 

\subsection{} 
The classification when $J=I$ (\ie in the case of
$\C[N]$ or $\C[B^-\backslash G]$) was given by Berenstein,
Fomin and Zelevinsky \cite{BFZ}. The only finite type
cases are $A_n\ (n\le 4)$.

\subsection{}
The classification when $J=\{j\}$ is a singleton 
and $G$ is of type $A_n$ (the Grassmannian $\Gr(j,n+1)$)
was already known \cite{S}.
Note that when $J=\{1\}$ (the projective space $\Pro^n(\C)$)
the cluster structure is trivial.
Indeed $\sub Q_{\{1\}}$ has only $n$ indecomposable objects
which are all $\Ext$-projective.
Thus ${\cal A}_{\{1\}}$ (resp. $\widetilde{\cal A}_{\{1\}}$)
has no cluster variable and is reduced to its coefficient
ring, a polynomial ring in $n$ (resp. $n+1$) variables.

\subsection{}
We now indicate how to obtain the classification.

\subsubsection{}\label{11.4.1}
The first step is to check that all cluster algebras 
${\cal A}_J$ of
Table~\ref{tab1} are indeed of finite type. 
For this, one chooses $\ii\in R(w_0,K)$ and computes following
\ref{ss84} the exchange matrix $B(\ii,J)$.
Let $\G_{\ii,J}$ denote the subgraph of $\G_\ii$ 
corresponding to the principal part of $B(\ii,J)$.
By \cite{FZ2}, one then has to find a sequence of mutations
transforming this graph into an orientation of the Dynkin
diagram of the claimed cluster type.
In each case, this is a straightforward verification.
We illustrate it in the next example.
\begin{example}\label{ex11.1}
{\rm
We continue Example~\ref{ex95}. Here $G$ is of type 
$A_5$ and $J=\{1,3\}$.
The graph $\G_{\ii,J}$ is shown in Figure~\ref{fig5}.
\begin{figure}[t]
\begin{center}
\leavevmode
\epsfxsize =6cm
\epsffile{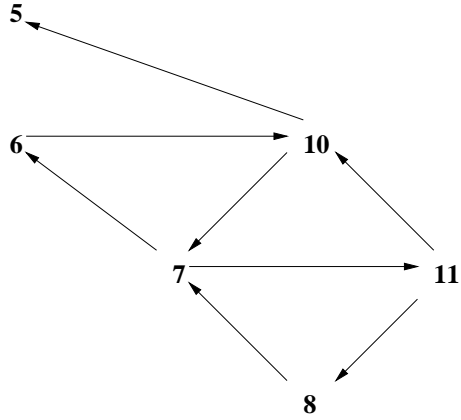}
\end{center}
\caption{\small {\it The graph $\Gamma_{\ii,J}$ for $J=\{1,3\}$ in type
$A_5$.}}
\label{fig5}
\end{figure}
If one performs a sequence of $3$ mutations, first at vertex $6$,
then at vertex $11$, finally at vertex $7$, one gets a quiver of
type $E_6$, in agreement with Table~\ref{tab1}.
\finex
}
\end{example}
For the case $D_n$, $J=\{n\}$, see also Section~\ref{sect12}
below.
For the case $D_5$, $J=\{1\}$, see also Section~\ref{sect13}
below.

\subsubsection{}\label{11.4.2}
One then needs to check a number of ``minimal'' infinite cases.
In type $A$ these are 
\begin{itemize}
\item $A_5\colon J=\{1,3,5\}$,
\item $A_6\colon J=\{1,3\}, \{1,4\}, \{3,4\}$,
\item $A_7\colon J=\{2,3\}, \{3,7\}$,
\item $A_n\ (n\ge 5)\colon J=\{2,n-1\}$.
\end{itemize}
In type $D$ these are
\begin{itemize}
\item $D_4\colon J=\{1,2,4\}$,
\item $D_5\colon J=\{1,2\}, \{1,5\}$,
\item $D_6\colon J=\{1\}$,
\item $D_n  (n\ge 4)\colon J=\{n-1\}$.
\end{itemize}
In type $E$, labelling the Dynkin diagrams as in \cite{Bo}, these are
\begin{itemize}
\item $E_6\colon J=\{1\}, \{2\}$.
\item $E_7\colon J=\{7\}$.
\item $E_8\colon J=\{8\}$.
\end{itemize}
These cases are settled by calculating as in \ref{11.4.1}
the graph $\G_{\ii,J}$. 
Then, one may either check that $\G_{\ii,J}$ contains a full
subgraph from the list of minimal infinite subgraphs
of Seven~\cite{Se}, or find a sequence of mutations transforming
$\G_{\ii,J}$ into a graph (containing a full subgraph) of
affine type. 
\begin{example}
{\rm 
Take $G$ of type $D_4$ and $J=\{3\}$.
Then 
\[
\ii=(1,2,4,3,1,2,4,3,1,2,4,3)
\]
belongs to $R(w_0,K)$.
The graph $\Gamma_{\ii,J}$ is shown in Figure~\ref{figD4}.
\begin{figure}[t]
\begin{center}
\leavevmode
\epsfxsize =4cm
\epsffile{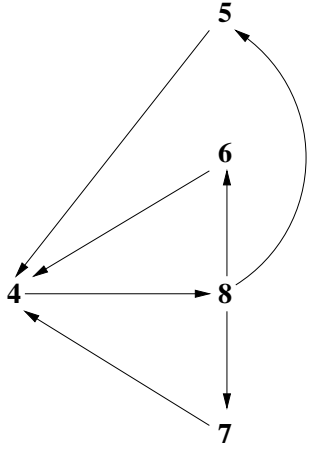}
\end{center}
\caption{\small {\it The graph $\Gamma_{\ii,J}$ for $J=\{3\}$ in type
$D_4$.}}
\label{figD4}
\end{figure}
If one performs a mutation at vertex~$4$,
one gets a quiver of
affine type $D_4$.
\finex
}
\end{example}
For the case $D_6$, $J=\{1\}$, see also Section~\ref{sect13}
below.

\subsubsection{}
To prove that there are no other finite type cluster algebras
${\cal A}_J$ than those listed in Table~\ref{tab1}, one uses
the following obvious 
\begin{Lem}\label{lem11.2}
\begin{itemize}
\item[(i)] Let ${\cal A}_J$ and ${\cal A}_{J'}$ be two algebras
attached to the same Dynkin diagram and suppose that $J\subset J'$.
If ${\cal A}_J$ has infinite cluster type then ${\cal A}_{J'}$
has infinite cluster type.
\item[(ii)] Let ${\cal A}_J$ and ${\cal A}'_J$ be two algebras
attached to two Dynkin diagrams $\De$ and $\De'$ and suppose
that $\De$ is a full subdiagram of $\De'$. 
If ${\cal A}_J$ has infinite cluster type then ${\cal A}'_J$
has infinite cluster type.
\end{itemize}\cqfd
\end{Lem}

\subsubsection{}
Assume that $G$ is of type $A_n$.
By \cite{S} we know that if $n\ge 8$
the only cluster algebras ${\cal A}_{\{j\}}$ of finite type are
obtained for $j=1,2,n-1,n$.
It follows from Lemma~\ref{lem11.2}~(i) that for $n\ge 8$ 
the algebra ${\cal A}_J$ can have
finite cluster type only if $J\subseteq\{1,2,n-1,n\}$.
Since $J=\{2,n-1\}$ yields an infinite type by \ref{11.4.2},
this is also the case for $J=\{1,2,n-1\}$ and $J=\{1,2,n-1,n\}$.
Hence, the only cases left for $n\ge 8$ are those of Table~\ref{tab1}.
For $A_7$, we need to check that no $J$ of the form $\{j,3\}$
gives a finite type.
By \ref{11.4.2}, $\{2,3\}$ and $\{3,7\}$ can be excluded.
Using Lemma~\ref{lem11.2}~(ii) we can also exclude $\{1,3\}$
and $\{2,3\}$ since they already give an infinite type for
$A_6$.
For $\{3,5\}$ and $\{3,6\}$, we use again Lemma~\ref{lem11.2}~(ii)
and restrict to type $A_6$ by removing the vertex $1$ of the
Dynkin diagram. This yields $J=\{2,4\}$ and $J=\{2,5\}$ in
type $A_6$. These can be ruled out by using \ref{11.4.2} again
($J=\{2,4\}$ gives already an infinite type for $A_5$).
This finishes type $A_7$.
The types $A_6$ and $A_5$ are dealt with similarly. 
Details are omitted. This finishes the classification in type $A$.

\subsubsection{}\label{tD}
Assume that $G$ is of type $D_n$.
If $n\ge 6$ and $j\not = n$ then ${\cal A}_{\{j\}}$ has infinite type.
This is easily shown by induction on $n$.
Indeed if $n=6$, by \ref{11.4.2} we can exclude $j=1,2,5$.
By Lemma~\ref{lem11.2}~(ii) and using again \ref{11.4.2} for $D_5$ 
and $D_4$ we can exclude $j=4$ and $j=3$. Thus $n=6$ is checked.
If the claim holds for $D_{n-1}$ then we can exclude 
$j=n-1$ by \ref{11.4.2} and use Lemma~\ref{lem11.2}~(ii) and
the assumption for $D_{n-1}$ for all other $j$'s.
If $n=5$, and ${\cal A}_{J}$ is of finite type, we must
have $J\subseteq \{1,2,5\}$.
By \ref{11.4.2} we can exclude all pairs in $\{1,2,5\}$, 
so we are left with $J=\{1\}$ or $J=\{5\}$ in agreement with
Table~\ref{tab1} ($J=\{2\}$ and $J=\{1\}$ are conjugate under
a diagram automorphism).
If $n=4$, we have a diagram automorphism of order $3$ exchanging
$1$, $2$ and $4$. Taking into account \ref{11.4.2}, we see that
we are left with the cases of Table~\ref{tab1} (up to isomorphism).
This finishes the classification in type $D$.

\subsubsection{}
Assume that $G$ is of type $E_6$.
We have a diagram automorphism exchanging $1$ and $6$.
Thus, using \ref{11.4.2}, we see that ${\cal A}_{\{j\}}$
is of infinite type for $j=1,2,6$. 
Using Lemma~\ref{lem11.2}~(ii) and \ref{tD}, we obtain
by reduction to $D_5$ that ${\cal A}_{\{j\}}$
is also of infinite type for $j=3,4,5$.
Thus ${\cal A}_{J}$ is infinite for every $J$.
The cases when $G$ is of type $E_7$ or $E_8$ follow
by means of Lemma~\ref{lem11.2}~(ii) and  \ref{11.4.2}.
This finishes the classification in type $E$.

\begin{figure}[t!]
\begin{center}
\leavevmode
\epsfxsize = 8.5cm
\epsffile{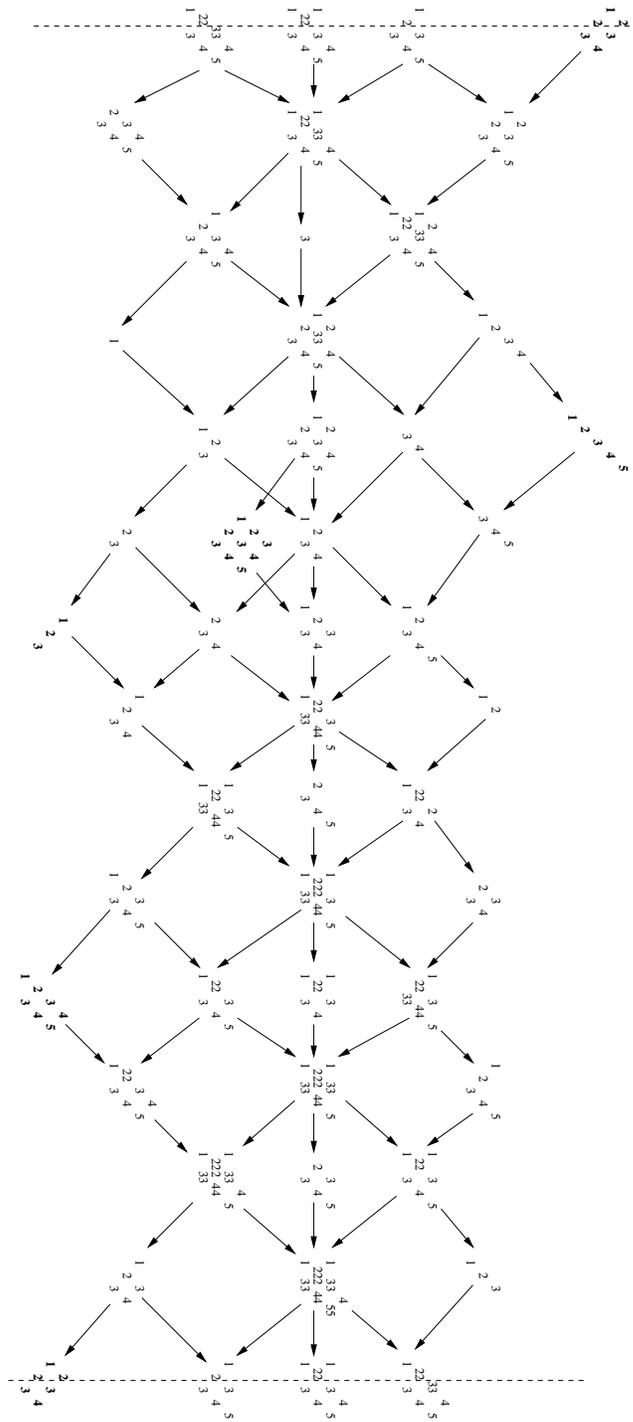}
\end{center}
\caption{\small {\it The Auslander-Reiten quiver for $J=\{1,3\}$ in type
$A_5$.}}
\label{figARA5J13}
\end{figure} 

\subsection{}
When ${\cal A}_{J}$ has infinite cluster type, it has 
an infinite number of cluster variables. 
Therefore the category $\sub Q_J$ has an infinite number of 
indecomposable rigid modules.

When ${\cal A}_{J}$ has finite cluster type, one can knit
the Auslander-Reiten quiver of $\sub Q_J$.
One obtains a finite connected quiver drawn either on a 
cylinder or on a M\"obius band. 
The corresponding stable Auslander-Reiten quiver, obtained
by deleting the Ext-projective modules, is isomorphic to
the Auslander-Reiten quiver of the cluster category of the same type
as ${\cal A}_{J}$, introduced by Buan, Marsh, Reineke, Reiten and 
Todorov \cite{BMRRT}.

\begin{example}
{\rm We continue Examples~\ref{ex95} and \ref{ex11.1}. 
Here $G$ is of type $A_5$ and $J=\{1,3\}$.
The Auslander-Reiten quiver of $\sub Q_J$ is shown in Figure~\ref{figARA5J13}.
This quiver is drawn on a M\"obius band. 
The horizontal dashed lines
at the top and at the bottom are to be identified after performing
a half-turn. The relative Auslander-Reiten translation $\tau_J$
is going up.  
Each indecomposable $\L$-module $M$ has a grading with semisimple 
homogeneous components. 
We represent $M$ by displaying the layers of this filtration.
Thus, 
\[
\begin{matrix}
&2\\1&&3\\&2&&4\\&&3
\end{matrix}
\]
stands for a $6$-dimensional module with graded components
$S_2$, $S_1\oplus S_3$, $S_2\oplus S_4$
and $S_3$, from top to bottom.
To determine this quiver, we have first calculated the
$\tau_J$-orbits using Proposition~\ref{prop2-2}, which states
that $\tau_J^{-1}$ is equal to the relative syzygy functor.
The five Ext-projective modules $L_i$
are printed in bold type. If one erases these five vertices as
well as the corresponding arrows, one obtains the stable
Auslander-Reiten quiver of $\sub Q_J$, which is isomorphic to
the quiver of a cluster category of type $E_6$.
\finex}
\end{example}

\section{The coordinate ring of a smooth
quadric in $\Pro^{2n-1}(\C)$}\label{sect12}

When $G$ is of type $D_n$ and $J=\{n\}$, the construction 
of the previous sections yields a cluster
algebra structure on the homogeneous coordinate ring of a smooth
quadric in $\Pro^{2n-1}(\C)$, and a finite type subcategory 
$\sub Q_n$ of $\mod \L$
which can be regarded as a lift of this ring.
We shall present this example in some detail. 

\subsection{}
Let $n\ge 4$.
Let $U=\C^{2n}$ and $\Pro=\Pro(U)=\Pro^{2n-1}(\C)$. 
Let $(u_1,\ldots,u_{2n})$ be a fixed basis of $U$
and $(y_1,\ldots,y_{2n})$ the coordinate functions with respect
to this basis.
We denote by $[y_1:y_2:\cdots : y_{2n}]$ the corresponding system 
of homogeneous coordinates on $\Pro$ and we consider the smooth 
quadric ${\mathcal Q}$ in $\Pro$ given by the equation 
\begin{equation}\label{eq6}
q(y_1,\ldots,y_{2n}) := \sum_{i=1}^n (-1)^i y_i\,y_{2n+1-i} = 0.
\end{equation}
In other words, ${\mathcal Q}$ is isomorphic to the variety
of isotropic lines in the quadratic space $(U,q)$.
Note that every smooth quadric in $\Pro$ can be brought to equation 
(\ref{eq6}) by an appropriate change of coordinates.
The homogeneous coordinate ring of ${\mathcal Q}$ is 
\[
\C[{\mathcal Q}]=\C[y_1,\ldots,y_{2n}]\slash(q(y_1,\ldots,y_{2n})= 0).
\]

\subsection{}
Let $G$ be the group of linear transformations of $U$ preserving
the quadratic form $q$. 
We regard $G$ as acting on the {\em right} on $U$,
\ie elements of $U$ are regarded as {\em row} vectors.
We identify $g\in G$ with its matrix with respect to 
$(u_1,\ldots,u_{2n})\colon$  the {\em i}th row of the matrix $g$
is the list of coordinates of $u_i g$.
The group $G=SO_{2n}(\C)$ is a group of type $D_n$.
In this realization, the subgroup of diagonal matrices of
$G$ is a maximal torus, and the subgroup of upper (\resp lower)
triangular matrices of $G$ is a Borel subgroup denoted
by $B$ (\resp $B^-$).
We label the vertices of the Dynkin diagram in such a way that
$1$ and $2$ correspond to the two spin representations and
$n$ to the vector representation $U$.

\subsection{}\label{ss10-3} 
Let $J=\{n\}$ and $K=\{1,2,\ldots,n-1\}$.
The quadric
${\mathcal Q}$ is isomorphic to $B_K^-\backslash G$. 
Indeed, $\C u_{1}$ is an isotropic line and ${\mathcal Q}$ is 
the $G$-orbit of $\C u_{1}$, the stabilizer of $\C u_{1}$ being
equal to $B_K^-$.
Therefore, for $g\in G$, the coordinates of the vector $u_1g$
give a system of homogeneous coordinates for the point 
$\C u_1g \in {\mathcal Q}$, that we can and shall continue
to denote by $y_1,\ldots ,y_{2n}$.
The affine open subset given by the non-vanishing of $y_1$ 
is isomorphic to $N_K$. Thus, setting $z_k = y_k/y_1$
\[
\C[N_K]=\C[z_2,\ldots,z_{2n}]\slash(q(1,z_2,\ldots,z_{2n})= 0).
\]

\subsection{}\label{ss10-4}
The functions on $G$ mapping $g$ to the coordinates of the vector
$u_1g$ are regular, and are nothing else than the generalized
minors $\De_{\varpi_n,v(\varpi_n)}$ for $v$ in the Weyl group.
In particular the coordinate of 
$u_{1}$ in $u_1g$ is equal to
$\De_{\varpi_n,\varpi_n}(g)$, and the coordinate of 
$u_{2n}$ in $u_1g$ is equal to
$\De_{\varpi_n,w_0(\varpi_n)}(g)$.
The restrictions $D_{\varpi_n,v(\varpi_n)}$ of all these functions
to the unipotent radical $N$ of $B$
are the elements $\varphi_{M}$ of $\C[N]$, where $M$ 
runs over the $2n$ submodules $M$ of the injective $\L$-module
$Q_n$.
In particular
\[
D_{\varpi_n,\varpi_n} = 1= \varphi_{\mathbf{0}} = z_1,\qquad
D_{\varpi_n,w_0(\varpi_n)} = \varphi_{Q_n} = z_{2n},
\]
where $\mathbf{0}$ denotes the trivial submodule of $Q_n$.

\subsection{} 
\begin{figure}[t]
\begin{center}
\leavevmode
\epsfxsize =9cm
\epsffile{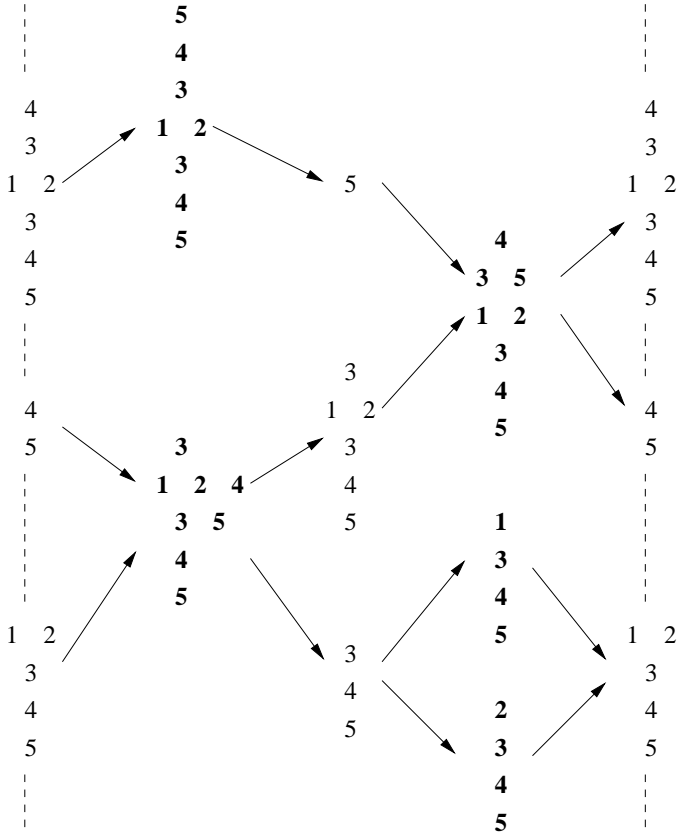}
\end{center}
\caption{\small 
{\it The Auslander-Reiten quiver of the category $\sub Q_5$ in type
  $D_5$.}}
\label{fig4}
\end{figure}

The category $\sub Q_n$ has a finite number of
indecomposable objects, which are all rigid. Here is a complete
list:
\begin{itemize}
\item the $2n-1$ nonzero submodules of $Q_n$, with pairwise distinct 
dimension vectors:
\begin{quote}
$
[0,\ldots,0,1],\ 
[0,\ldots,0,1,1],\ 
\ldots,\
[0,1,\ldots,1,1],\ 
$

$
[1,0,1,\ldots,1,1],\
[1,1,1,\ldots,1,1],\
[1,1,2,1,\ldots,1,1],\
$

$
[1,1,2,2,1,\ldots,1,1],
\ldots,\,
[1,1,2,2,\ldots,2].
$
\end{quote}
The modules with dimension vectors 
$[0,1,\ldots,1,1],\ [1,0,1,\ldots,1,1], [1,1,2,2,\ldots,2]$
are the relatively projective modules $L_2, L_1, L_n=Q_n$,
respectively.
We shall denote the regular functions on $N$ corresponding 
to these $2n-1$ modules by 
$z_2, \ldots, z_{2n}$, in agreement with \ref{ss10-3} and \ref{ss10-4}. 
They generate $\C[N_K]$ and satisfy $q(1,z_2,\ldots,z_{2n})=0$.

\item the $n-3$ indecomposable submodules of $Q_n\oplus Q_n$
with socle $S_n\oplus S_n$ (up to isomorphism). They all have
the same dimension vector as $Q_n$, namely $[1,1,2,2,\ldots,2]$,
and they all are projective objects in $\sub Q_n$.
They form the remaining indecomposable projectives 
$L_3, \ldots , L_{n-1}$ of $\sub Q_n$.
We denote the corresponding functions on $N$ by $p_3, \ldots,
p_{n-1}$.
One can check that for $3\le k\le n-1$,
\[
p_k=z_{n+1-k}z_{n+k} - z_{n-k}z_{n+k+1} + \cdots 
+ (-1)^{n-k-1}z_2z_{2n-1} + (-1)^{n-k}z_{2n}.
\]
One can also check that $p_k=D_{u_k(\varpi_k),w_0(\varpi_k)}$
for some appropriate element $u_k$ of $W$.
\end{itemize}
As an illustration,
the Auslander-Reiten quiver of the category $\sub Q_5$ in type $D_5$ is
displayed in Figure~\ref{fig4}.
The quiver is drawn on a cylinder, obtained by identifying the two
vertical dashed lines.
The 5 projective objects are written in bold type. 

\subsection{} The cluster algebra structure ${\cal A}_J$ on 
$\C[N_K]$ is obtained as follows. 
The coefficient ring is generated by $p_1, \ldots,p_n$,
where $p_1:=z_n,\ p_2:=z_{n+1},\ p_n:=z_{2n}$ and the remaining
$p_k$'s are defined above.
There are $2(n-2)$ cluster variables, namely
\[
z_2,\ z_3, \ldots,\ z_{n-1},\ z_{n+2},\ldots ,z_{2n-1}.
\]
The first $n-2$ variables form a cluster and for $2\le k\le n-1$
we have the exchange relations
\[
z_kz_{2n-k+1} = 
\left\{
\begin{array}{ll}
p_{n-k+1}+p_{n-k+2} & \mbox{if $2\le k \le n-2$,}\\[2mm]
p_1p_2+p_3 & \mbox{if $k=n-1$.}
\end{array}
\right.
\]
This shows that ${\cal A}_J=\C[N_K]$ has finite cluster type 
equal to $(A_1)^{n-2}$,
in agreement with Table~\ref{tab1}.

\subsection{}
The cluster algebra structure $\widetilde{\cal A}_J$ on 
$\C[{\mathcal Q}]$ is obtained as follows. 
The coefficient ring is generated by $q_0,q_1, \ldots,q_n$,
where $q_0:=y_1,\ q_1:=y_n,\ q_2:=y_{n+1},\ q_n:=y_{2n}$ and 
for $3\le k\le n-1$,
\[
q_k:=y_{n+1-k}y_{n+k} - y_{n-k}y_{n+k+1} + \cdots 
+ (-1)^{n-k-1}y_2y_{2n-1} + (-1)^{n-k}y_1y_{2n}.
\]
There are $2(n-2)$ cluster variables, namely
\[
y_2,\ y_3, \ldots,\ y_{n-1},\ y_{n+2},\ldots ,y_{2n-1}.
\]
The first $n-2$ form a cluster and for $2\le k\le n-1$
we have the exchange relations
\[
y_ky_{2n-k+1} = 
\left\{
\begin{array}{ll}
q_{n-1}+q_0q_n & \mbox{if $k=2$,}\\[2mm]
q_{n-k+1}+q_{n-k+2} & \mbox{if $3\le k \le n-2$,}\\[2mm]
q_1q_2+q_3 & \mbox{if $k=n-1$.}
\end{array}
\right.
\]
Thus, $\widetilde{\cal A}_J = \C[{\mathcal Q}]$
 is also a cluster algebra of type $(A_1)^{n-2}$.

\subsection{}
When $n=3$, $\C[{\mathcal Q}]$
has a cluster algebra structure of type $A_1$.
Using the same notation as above for the generators
$y_1, \ldots, y_6$ and the coefficients $q_0, \ldots, q_3$,
the unique exchange relation reads 
\[
y_2y_5=q_1q_2+q_0q_3.
\]
Note that $D_3 \cong A_3$ and that ${\mathcal Q}$ is
isomorphic to the Grassmannian of $2$-planes in $\C^4$.

\section{Isotropic Grassmannians}\label{sect13}

We retain the notation of Section~\ref{sect12}.
Let ${\cal G}$ denote the Grassmann variety of totally
isotropic $n$-subspaces of $U$.
This variety has two connected components, and we shall
denote by ${\cal G}_0$ the component containing the subspace
spanned by $(u_1,\ldots,u_n)$.
The stabilizer of this subspace for the natural action
of $G$ is the maximal parabolic subgroup $B_K^-$,
where now $K=\{2,\ldots ,n\}$, that is, $J=\{1\}$.
Thus ${\cal G}_0$ is isomorphic to $B_K^-\backslash G$.
We shall now discuss the cluster algebra structure on 
$\C[{\cal G}_0]$ and the corresponding subcategory 
$\sub Q_1$.

\subsection{}\label{ss13.3}
\begin{figure}[t]
\begin{center}
\leavevmode
\epsfxsize =10cm
\epsffile{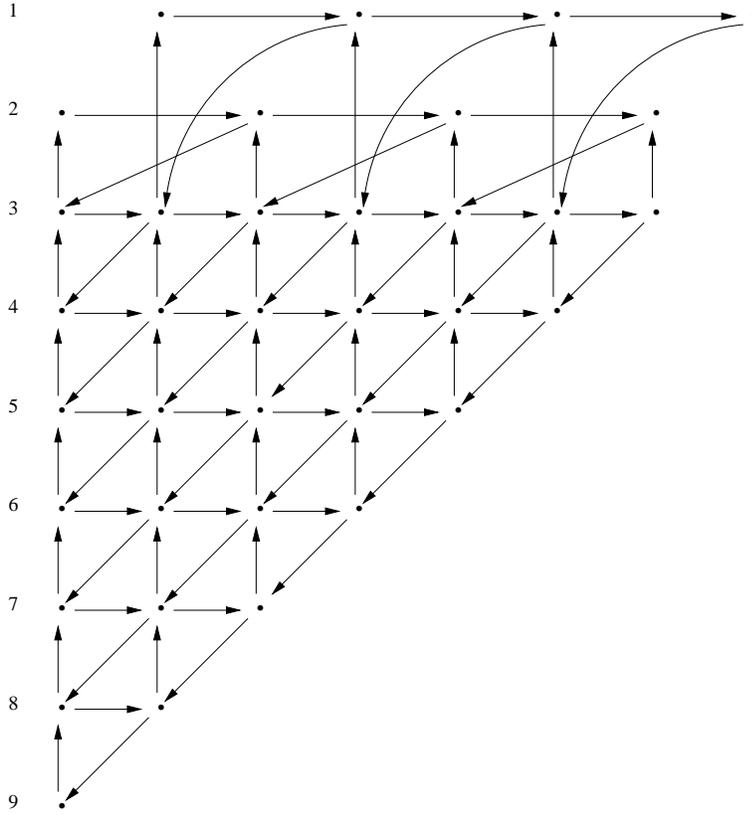}
\end{center}
\caption{\small 
{\it The triangular grid for $D_9$.}}
\label{tgrid}
\end{figure}

As in \ref{grid}, the rule of \ref{ss84} gives us a unique 
initial cluster for ${\cal A}_{\{1\}}$. 
The graph encoding the exchange matrix of this cluster
has the shape of a triangular grid in which
the first two rows, corresponding to the exceptional vertices 1 and 2 of the
Dynkin diagram of $D_n$, have a special structure.
A typical example ($n=9$) is displayed in Figure~\ref{tgrid}
(compare Figure~\ref{fig2}, the rectangular grid attached
to the ordinary Grassmannian).
The $n$ generators of the coefficient ring correspond to the
leftmost vertices on each row.
It is easy to see that this graph yields a cluster algebra of
infinite type, except for $n=4$ and $n=5$.

\subsection{}
When $n=4$, because of the order 3 symmetry of the Dynkin diagram,
${\cal G}_0$ is isomorphic to the quadric ${\cal Q}$ of
Section~\ref{sect12}, and $\sub Q_1$ is equivalent to $\sub Q_4$.
In particular $\widetilde{\cal A}_{\{1\}}=\C[{\cal G}_0]$ is a 
cluster algebra of type $A_1\times A_1$.

\subsection{}
\begin{figure}[t]
\begin{center}
\leavevmode
\epsfxsize =14cm
\epsffile{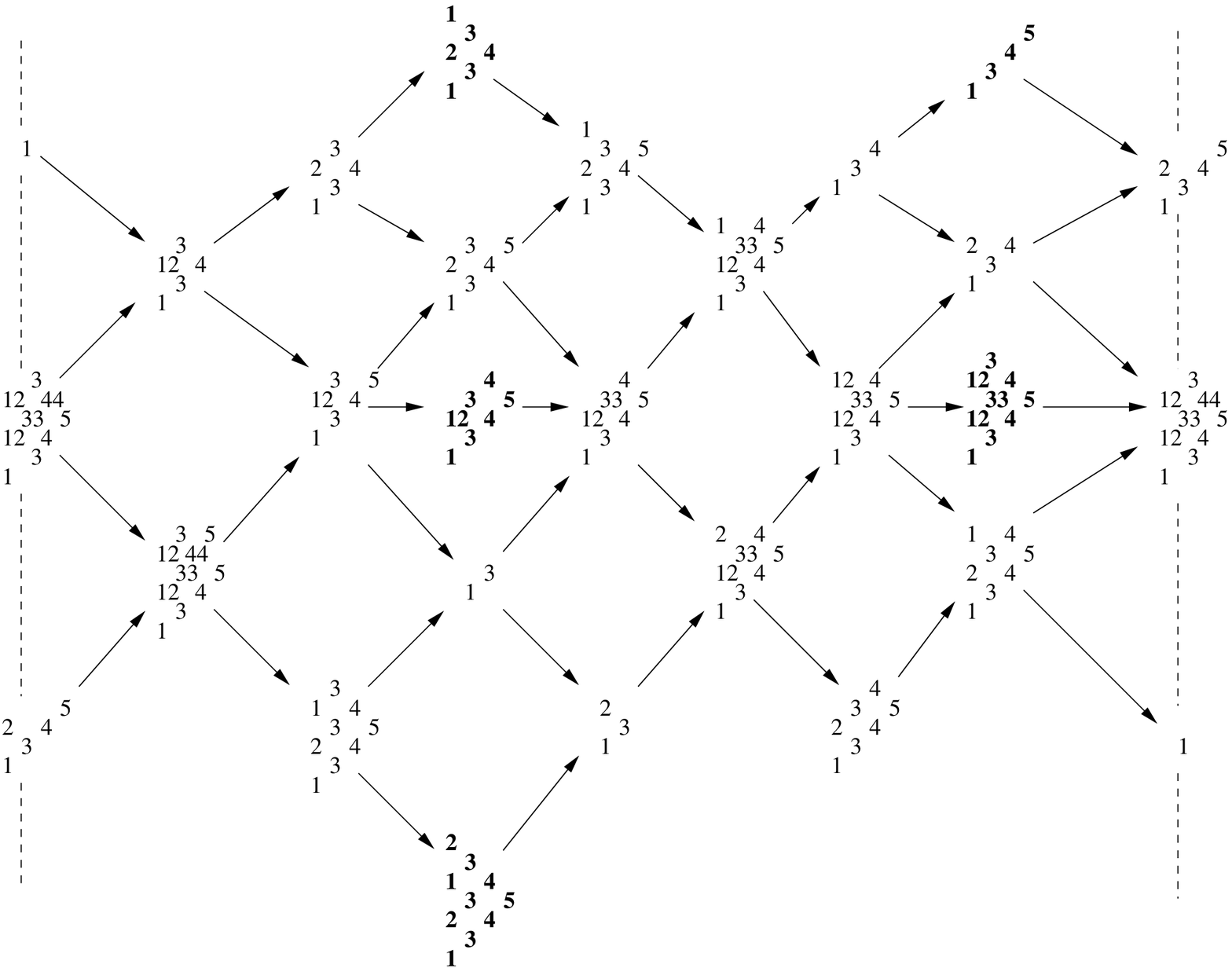}
\end{center}
\caption{\small 
{\it The Auslander-Reiten quiver of the category $\sub Q_1$ in type
  $D_5$.}}
\label{ARSubQ1D5}
\end{figure}

When $n=5$, $\sub Q_1$ is a category of finite type with $25$
indecomposable objects (up to isomorphism), $5$ of them being 
$\Ext$-projective.
The Auslander-Reiten quiver of $\sub Q_1$ is displayed in
Figure~\ref{ARSubQ1D5}. 
It is drawn on a M\"obius band obtained by identifying the 
two vertical dashed lines after performing a half-turn.
The stable Auslander-Reiten quiver (obtained by deleting
the $5$ projectives) is the quiver of a cluster
category of type $A_5$, in agreement with Table~\ref{tab1}.

We shall now describe the cluster algebra structure 
$\widetilde{\cal A}_{\{1\}}$ on $\C[{\cal G}_0]$.
First we express the cluster 
variables and the generators of the coefficient ring in terms of 
the generalized minors $\De_{\varpi_1,u(\varpi_1)}$
which generate $\C[{\cal G}_0]$.
(Note that, strictly speaking, these minors are not functions on 
$G=SO_{2n}(\C)$ but only on the corresponding simply connected
group ${\rm Spin}_{2n}(\C)$.)
The $16$ minors $\De_{\varpi_1,u(\varpi_1)}$ will be denoted
for short by $\De_i\ (1\le i\le 16)$
according to Table~\ref{tableminD5}.
\begin{table}
\begin{center}
\begin{tabular}
{|c|c|c|}
\hline
$i$& $u$ & mutation sequence for $\De_i$\\
\hline
1&$e$ & $q_0$\\
2&$s_1$  & $z_1$\\ 
3&$s_3s_1$  & $z_2$\\
4&$s_2s_3s_1$  & $z_3$\\
5&$s_4s_3s_1$ &  $z_4$\\
6&$s_2s_4s_3s_1$ & $\mu_2$\\
7&$s_5s_4s_3s_1$ & $q_5$\\
8&$s_3s_2s_4s_3s_1$ & $\mu_1$\\
9&$s_2s_5s_4s_3s_1$ & $\mu_2\mu_1\mu_5\mu_4$\\
10&$s_1s_3s_2s_4s_3s_1$ & $q_1$\\
11&$s_5s_3s_2s_4s_3s_1$ & $\mu_1\mu_5\mu_4$\\
12&$s_1s_5s_3s_2s_4s_3s_1$ & $\mu_5\mu_4$\\
13&$s_4s_5s_3s_2s_4s_3s_1$ & $\mu_4\mu_3\mu_2\mu_1\mu_5\mu_4$\\
14&$s_1s_4s_5s_3s_2s_4s_3s_1$ & $\mu_1\mu_4\mu_3\mu_2\mu_1\mu_5\mu_4$\\
15&$s_3s_1s_4s_5s_3s_2s_4s_3s_1$ & $\mu_3\mu_1\mu_5\mu_4$\\
16&$s_2s_3s_1s_4s_5s_3s_2s_4s_3s_1$ & $q_2$\\
\hline
\end{tabular}
\end{center}
\caption{\small \it Generalized minors $\De_{\varpi_1,u(\varpi_1)}$
in type $D_5$.
\label{tableminD5}} 
\end{table} 
Among them, 
\[
\De_1=q_0,\quad 
\De_7=\widetilde{\varphi_{L_5}}=q_5,\quad
\De_{10}=\widetilde{\varphi_{L_1}}=q_1,\quad
\De_{16}=\widetilde{\varphi_{L_2}}=\widetilde{\varphi_{Q_1}}=q_2,
\]
are generators of the coefficient ring of $\widetilde{\cal A}_{\{1\}}$.
The two other generators of the coefficent ring are
\[
q_3=\widetilde{\varphi_{L_3}}=\De_4\De_{15}-\De_3\De_{16}, \qquad
q_4=\widetilde{\varphi_{L_4}}=\De_2\De_{14}-\De_1\De_{13}.
\]
The initial cluster, obtained by lifting to 
$\widetilde{\cal A}_{\{1\}}$
the cluster of \ref{ss13.3},  
consists of the functions
\[
z_1=\De_2,\quad
z_2=\De_3,\quad
z_3=\De_4,\quad
z_4=\De_5,\quad
z_5=\De_2\De_8-\De_1\De_{10}.
\]
The exchange matrix of this cluster is
\[
\tB=
\begin{pmatrix}
         0&0&0&0&1\cr
         0&0&1&1&-1\cr
         0&-1&0&0&0\cr
         0&-1&0&0&1\cr
         -1&1&0&-1&0\cr
           \hline        
         1&0&0&0&-1\cr
         0&0&1&0&0\cr
         0&0&-1&0&1\cr
         0&0&0&1&-1\cr
         0&0&0&-1&0\cr
         \hline
         1&0&0&0&0
\end{pmatrix},
\]
where the successive rows are labelled by 
$z_1, z_2, z_3, z_4, z_5, q_1, q_2, q_3, q_4, q_5, q_0$.
The last column of Table~\ref{tableminD5} indicates which 
sequence of mutations produces, starting from this initial cluster,
each minor $\De_i$. 
This shows that $\widetilde{\cal A}_{\{1\}}$ 
contains a set of generators of $\C[{\cal G}_0]$.
Hence, Conjecture~\ref{conj1} and Conjecture~\ref{conj2} 
are also proved in this case.
(Note that $\varpi_1$ is a minuscule weight, so 
no localization is needed in Conjecture~\ref{conj2}.)

The remaining cluster variables all have degree $2\varpi_1$ and
are given by the following quadratic expressions in the minors $\De_i$:
\[
\De_3\De_{13}-\De_1\De_{15},\quad
\De_4\De_{14}-\De_2\De_{16},\quad
\De_6\De_{15}-\De_5\De_{16},\quad
\De_2\De_{11}-\De_1\De_{12},\quad
\]
\[
\De_4\De_{13}-\De_1\De_{16},\quad
\De_9\De_{15}-\De_7\De_{16},\quad
\De_5\De_{12}-\De_7\De_{10}.
\]
As already mentioned, $\widetilde{\cal A}_{\{1\}}=\C[{\cal G}_0]$ 
is a cluster algebra of finite type $A_5$.

\subsection{}
\begin{figure}[t]
\begin{center}
\leavevmode
\epsfxsize =9cm
\epsffile{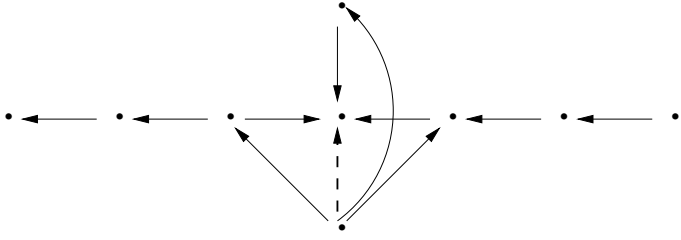}
\end{center}
\caption{\small 
{\it The elliptic diagram of type $E_7^{(1,1)}$.}}
\label{E_7^11}
\end{figure}

When $n=6$, there exists a sequence of 6 mutations transforming
the principal part of the triangular grid into the graph 
displayed in Figure~\ref{E_7^11}.
Here the dashed arrow stands for a pair of entries $\pm 2$
in the exchange matrix. 
Following \cite{GLS0}, we propose to attach to this 
infinite type cluster algebra $\widetilde{\cal A}_{\{1\}}$
the elliptic type $E_7^{(1,1)}$.

\section{Remarks on the non simply-laced case}

\subsection{}
Although the definition of the cluster algebra ${\cal A}_J$
was obtained using the representation theory of
preprojective algebras, the definition of the initial
seed in \ref{ss84} can be formulated without any reference
to preprojective algebras.
As a result, the same definition can serve to introduce
similar cluster algebras ${\cal A}_J$ in the 
non simply-laced types.
(This was suggested to us by Andrei Zelevinsky.)
One can expect that ${\cal A}_J$ is again equal
to $\C[N_K]$ and can be lifted to a cluster algebra
structure $\widetilde{\cal A}_J$ on $\C[B^-_K\backslash G]$,
where $G$ is now the corresponding algebraic group of
non simply-laced type.
 
\subsection{}
A similar study as in Section~\ref{sect11} gives 
the classification of all finite type cluster algebras
${\cal A}_J$ in the non simply-laced case.
The results are summarized in Table~\ref{tab2}.
Our convention for labelling the Dynkin diagrams of type
$B_n$ and $C_n$ is that the vertex associated with the
vector representation is numbered $n$. 
\begin{table}[h]
\begin{center}
\begin{tabular}
{|c|c|c|}
\hline
Type of $G$ & $J$ & Type of ${\cal A}_J$\\
\hline
$B_n$ $(n\ge 2)$& $\{n\}$ & $(A_1)^{n-1}$\\
\hline
$C_n$ $(n\ge 2)$& $\{n\}$ & $(A_1)^{n-1}$ \\
\hline
$B_2=C_2$& $\{1,2\}$ & $B_2=C_2$\\
\hline
$B_3$& $\{1\}$ & $C_3$\\
\hline
$C_3$& $\{1\}$ & $B_3$\\
\hline
\end{tabular}
\end{center}
\caption{\small \it Algebras ${\cal A}_J$ of finite cluster type (non
  simply-laced case).
\label{tab2}}
\end{table}


\bigskip\noindent
{\bf Acknowledgements.\ }
This paper was written during a stay at the Mathematisches
Forschungs\-insti\-tut Oberwolfach in July-August 2006.
We are very greatful to this institution for its support,
its hospitality and for providing ideal working conditions.
Some preliminary work was done while B.~Leclerc
was participating in the program {\em Algebraic Combinatorics}
at the Mittag-Leffler Institute (Stockholm, April 2005),
and in the program {\em Group Representation Theory}
at the Centre Bernoulli (Lausanne, June 2005).
He wants to thank the organizers of these programs for
inviting him.


\bigskip
\small

\noindent
\begin{tabular}{ll}
Christof {\sc Gei{ss}} : &
Instituto de Matem\'aticas,\\ 
&Universidad Nacional Aut\'onoma de M\'exico\\
& 04510 M\'exico D.F., M\'exico \\
&email : {\tt christof@math.unam.mx}\\[5mm]
Bernard {\sc Leclerc} :&
Universit\'e de Caen, LMNO UMR 6139\\
& 14032 Caen cedex, France\\
&email : {\tt leclerc@math.unicaen.fr}\\[5mm]
Jan {\sc Schr\"oer} :&
Mathematisches Institut, Universit\"at Bonn,\\
&Beringstr. 1, D-53115 Bonn, Germany\\
&email : {\tt schroer@math.uni-bonn.de}
\end{tabular}

\end{document}